\documentclass[a4paper,12pt]{article}
\usepackage[margin=1.2in]{geometry}
\pdfoutput=1
\usepackage{mathtools}

\usepackage[utf8]{inputenc}
\usepackage{amsfonts,amssymb, bm, amsthm}
\usepackage{graphicx,color,epstopdf}
\allowdisplaybreaks
\newtheorem{mylemma}{Lemma}[section]
\newtheorem{mytheorem}{Theorem}[section]
\newtheorem{mycorollary}{Corollary}[section]
\newtheorem{myremark}{Remark}[section]
\newtheorem{mydef}{Definition}[section]

\def\XXint#1#2#3{{\setbox0=\hbox{$#1{#2#3}{\int}$}
    \vcenter{\hbox{$#2#3$}}\kern-.5\wd0}}

\def\hG {\hat{G}}
\def\tbeta {\tilde{\beta}}
\def\tx {\tilde{x}}
\def\tX {\tilde{X}}
\def\tY {\tilde{Y}}
\def\ty {\tilde{y}}

\def\tu {\tilde{u}}
\def\tU {\tilde{U}}

\def\dd {''}

\def\hG {\hat{G}}

\def\bi{{\bf i}}

\newcommand{\tcr}{\textcolor{black}}

\begin{document}
\title{Wave scattering in layered orthotropic media I: a stable PML and a
  high-accuracy boundary integral equation method} \author{Yong Gao$^1$ and
  Wangtao Lu$^2$} \footnotetext[1]{School of Mathematical Sciences, Zhejiang
  University, Hangzhou 310027, China. Email: xzgaoyong@126.com.}
\footnotetext[2]{Corresponding author. School of Mathematical Sciences, Zhejiang
  University, Hangzhou 310027, China. Email: wangtaolu@zju.edu.cn. This author
  is partially supported by NSFC Grant 12174310 and by NSF of Zhejiang Province
  for Distinguished Young Scholars (LR21A010001). }

\maketitle
\begin{abstract}
  In anisotropic media, the standard perfectly matched layer (PML) technique suffers irrevocable instability in terminating the unbounded problem domains. It remains an open question whether a stable PML-like absorbing boundary condition exists. For wave scattering in a layered orthotropic medium, this question is affirmatively answered for the first time in this paper. In each orthotropic medium, the permittivity tensor uniquely determines a change of coordinates, that transforms the governing anisotropic Helmholtz equation into an isotropic Helmholtz equation in the new coordinate system. This leads us to propose a novel Sommerfeld radiation condition (SRC) to rigorously characterize outgoing waves in the layered orthotropic medium. Naturally, the SRC motivates a regionalized PML (RPML) to truncate the scattering problem, in the sense that a standard PML is set up in the new coordinate system in each orthotropic region. It is revealed that the RPML is unconditionally stable compared with the unstable uniaxial PML. A high-accuracy boundary integral equation (BIE) method is developed to solve the resulting boundary value problem. Numerical experiments are carried out to validate the stability of the RPML and the accuracy of the BIE method, showing exponentially decaying truncation errors as the RPML parameters increase.
\end{abstract}

\section{Introduction}
Wave scattering problems arise from a wide range of realistic applications
\cite{chew95} including optics, radar, remote sensing, seismology, etc.
Efficient and accurate numerical methods are highly desirable in related
simulations. One essential difficulty among many others is how to accurately
truncate the unbounded problem domains in the first place. Coined and proposed
by Berenger \cite{ber94} in 1994, the perfectly matched layer (PML) technique
has since then become a widespread truncation approach, due to its nearly zero
reflection, easy implementation, and friendly incorporation into standard
numerical methods \cite{tofhug00,mon03}. Originally, PML was used to terminate
waves in homogeneous and isotropic background media, and subsequently, has so
far been successfully extended to more complicated backgrounds, such as layered
\cite{luluqia18,luluson19} and periodic \cite{chewu03,yuhulurat21} structures,
still isotropic though.

In contrast, significant difficulty arises in the extension of PML to
anisotropic backgrounds, as reported in the literature
\cite{becfaujol03,skeadacra07, ton15, bonfliton18}. Roughly speaking, PML is
only able to absorb waves with outgoing group velocities, but cannot work for
anisotropic media with two group velocities of different signs, e.g., in
elastodynamics, since waves can even exponentially blow up in a PML region.
B\'ecache et al. \cite{becfaujol03} rigorously studied the instability of PML in
the time domain of elastic waves in anisotropic media, and derived a high
frequency stability (HFS) condition of the PML. Even worse, they further pointed
out that PML for electro-magnetic (EM) or acoustic waves in orthotropic media is
unconditionally instable. To circumvent PML, Bonnet-BenDhia et al.
\cite{bonfliton18} developed a novel half-space matching method to solve an
acoustic wave scattering problem in an orthotropic medium with a bounded
obstacle. By using analytic integral representations in half spaces surrounding
the obstacle, they found an alternative artificial boundary condition on the
edges of the half spaces to truncate the unbounded domain, which have great
potentials in handling more general scattering problems
\cite{bonchaflihazpertja21}. Nevertheless, this boundary condition is, arguably,
not quite easy to implement compared with PML. Consequently, an attractive but
still open question is: Can we find a stable PML-like absorbing boundary
condition in anisotropic media? \tcr{ It is worthwhile to mention the excellent
  work of Demaldent and Imperiale \cite{demimp13}, where they designed a stable
  PML through the use of changes of coordinates to solve an acoustic scattering
  problem with a homogeneous and anisotropic background, the same structure as
  in \cite{bonfliton18}. Nevertheless, their principle of designing stable PMLs
  relies on the aforementioned HFS condition for homogeneous anisotropic
  backgrounds in \cite{becfaujol03}. To the authors' best knowledge, an HFS
  condition for layered anisotropic media is still absent, so that it remains
  open to design a stable PML and to develop a theory to justify the stability
  of a PML for such layered backgrounds. Motivated by this, this paper studies a
  TM-polarized EM wave of frequency domain propagating in a two-dimensional (2D)
  two-layer orthotropic medium, and designs a regionalized PML (RPML), the
  stability of which can be verified easily with the aid of newly defined
  outgoing waves.}

To design a PML-like absorbing layer, the primary task is to understand the
radiation behavior of the total wave field, the nonzero component of the
magnetic field, at infinity. To achieve this, in either orthotropic layer, we
elaborately choose a special change of coordinates to transform the 2D governing
equation, an anisotropic Helmholtz equation, into an isotropic Helmholtz
equation, so that the interface between the two layers are perfectly matched.
Unlike the usual transition matrix for orthotropic media \cite{colkremon97}, our
transition matrix differs by an orthogonal matrix multiplier, which uniquely
depends on the anisotropicities of the two involving media.  The specially designed changes of coordinates
make the method of Fourier transform applicable to determine the background
Green's function. By studying the radiation behavior of the Green's function at
infinity and by establishing the fundamental Green's representation formula, we
are able to propose a novel Sommerfeld radiation condition (SRC) for the total
wave field, to define what an outgoing wave is in orthotropic media rigorously
\tcr{(see Definition 2.1)}, and to determine the far-field pattern of any
outgoing wave conveniently. Certainly, the Green's representation formula
directly induces an exact transparent boundary condition to truncate the
unbounded domain, but rapidly evaluating the costly background Green's function
and its derivatives requires carefully designed algorithms \cite{cai13}.

Instead, the SRC motivates us to design an RPML to truncate the unbounded
domain, in the sense that a standard PML is set up regionally in the new
coordinate system in either orthotropic layer. The new RPML inherits all the
advantages of the classical PML. \tcr{More importantly, the outgoing behavior of
  the Green's function and the related Green's representation formula imply that
  our RPML is unconditionally stable to perfectly absorb any outgoing waves
  defined in Definition 2.1.} Consequently, all standard numerical methods can now be
incorporated readily. To design a high-accuracy solver, we adopt our previously
developed PML-based boundary integral equation (BIE) method \cite{luluqia18}.
Numerical experiments are carried out to validate the stability of the RPML and
the high-accuracy of the BIE method. Numerical results show that the truncation
error due to the RPML decays exponentially as the RPML parameters increase.


The rest of this paper is organized as follows. In Section 2, we present a
mathematical formulation and propose a radiation condition for the scattering
problem. In Section 3, we apply the method of Fourier transform to compute the
background Green's function and study its asymptotic behavior at infinity. In
Section 4, we derive the Green's representation formula and the far-field
pattern of an outgoing wave. In Section 5, we present the setup of an RPML and
the implementation of a high-accuracy BIE method, and analyze the stability of
the RPML. In Section 6, we study several numerical examples. In Section 7, we
conclude this paper and present some potential extensions of the RPML technique.
\section{Problem formulation}
\begin{figure}[!htb]
\centering
\vspace{-0.4cm}
\includegraphics[width=6cm]{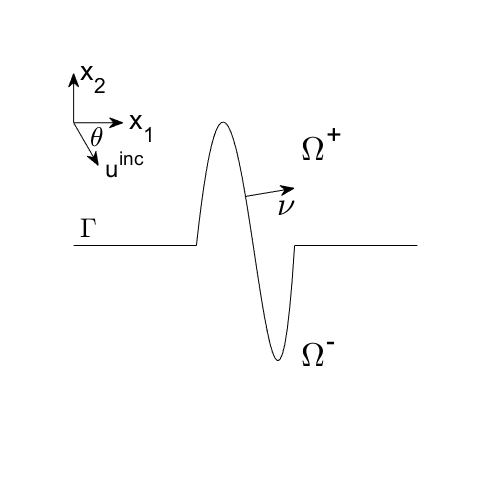}
\vspace{-1.9cm}
\caption{Profile of a 2D two-layer orthotropic medium. A locally perturbed
  straight line $\Gamma$ separate two unbounded regions $\Omega^+$ and
  $\Omega^-$, $\theta$ indicates the incident angle of a plane incident wave
  $u^{\rm inc}$, and $\nu$ is the unit normal vector of $\Gamma$ towards
  $\Omega^+$.}
\label{fig:layer}
\end{figure}

Let $(x_1,x_2,x_3)^{T}$ denote the Cartesian coordinate system of the
three-dimensional space $\mathbb{R}^3$. As shown in Figure~\ref{fig:layer},
two homogeneous and orthotropic media, which are invariant in the
$x_3$-direction, occupy two locally perturbed half planes $\Omega^+$ and
$\Omega^-$, respectively, in the $x_1Ox_2$-plane. The permittivities in
$\Omega^{\pm}$ are characterized by the following two positive definite tensors
\begin{equation*}
  \epsilon^\pm=
  \begin{bmatrix}

     \epsilon_{11}^\pm & \epsilon_{12}^\pm & 0 \\

     \epsilon_{21}^\pm & \epsilon_{22}^\pm & 0 \\
     
     0 & 0 & \epsilon_{33}^\pm \\
	\end{bmatrix},
\end{equation*}
respectively. In this paper, we assume that the locally perturbed straight line
$\Gamma$ separating $\Omega^\pm$ is on $x_2=0$ and is piecewise smooth.
\tcr{Throughout this paper, we assume that at least one of the two matrices is
non-diagonal, i.e. $|\epsilon^+_{12}| + |\epsilon^-_{12}|\neq 0$, to avoid trivial situations. }

A generic time-harmonic EM field propagating in
$\Omega^{\pm}$ is governed by Maxwell's equations:
\begin{align}
  \label{eq:max:1}
  \nabla_3 \times {\bm E} -\bi k_0 {\bm H}&=0, \\
  \label{eq:max:2}
  \nabla_3 \times {\bm H} +\bi k_0 \epsilon(x) {\bm E}&=0,
\end{align}
where $\nabla_3=(\partial_{x_1},\partial_{x_2},\partial_{x_3})^{T}$, $x=(x_1,x_2)^{T}$, ${\bm E}$
is the total electric field, ${\bm H}$ is the (properly scaled) total magnetic
field, $k_0=\frac{2\pi}{\lambda}$ is the free-space wavenumber, $\lambda$ is
the free-space wavelength, and
\begin{align}
  \epsilon(x) := [\epsilon_{ij}(x)]_{3\times 3}=\begin{cases}
    \epsilon^+, & x\in\Omega^+,\\
    \epsilon^-, & x\in\Omega^-.
    \end{cases}
\end{align}
Such vectorial equations can be further simplified if the EM field possesses
certain polarization. In this paper, we consider TM polarization only, i.e.,
${\bm H}= (0,0,u^{\rm tot})^{T}$ where the nonzero function $u^{\rm tot}$ is
assumed to be $x_3$-independent. Since
\begin{equation*}
  \nabla_3\times [\epsilon(x)^{-1}\nabla_3\times {\bm H}] - k_0^2 {\bm H} = 0,
\end{equation*}
$u^{\rm tot}$ satisfies the following anisotropic Helmholtz equation,
\begin{equation}
  \label{eq:Hhel}
\nabla \cdot (M(x)\nabla u^{\rm tot}(x))+k_0^{2}u^{\rm tot}(x)=0,
\end{equation}
for $x\in\mathbb{R}^2\backslash\Gamma$, where $\nabla =
(\partial_{x_1},\partial_{x_2})^{T}$,
\begin{equation*}	
M(x)=\frac{1}{\epsilon_{11}(x)\epsilon_{22}(x)-\epsilon_{12}(x)^2}\begin{bmatrix}

     \epsilon_{11}(x) & \epsilon_{12}(x)  \\

     \epsilon_{12}(x) & \epsilon_{22}(x) 
	\end{bmatrix}
\end{equation*}
is positive definite for any $x\notin \Gamma$. Across $\Gamma$, the following
continuous condition
\begin{align}
  \label{eq:cont:cond}
  [u^{\rm tot}]_{\rm j}= [\nu\cdot M\nabla u^{\rm tot}]_{\rm j}=0,
\end{align}
holds, where $\nu$ denotes the unit normal vector of $\Gamma$ towards $\Omega^+$,
and $[\cdot]_{\rm j}$ indicates the jump of the quantity.

For $x\in\Omega^{\pm}$, let $M_\pm=M(x)$,
\begin{equation}
  \label{eq:Mpm}
  M_{\pm}^{-1/2} = \left[
    \begin{array}{cc}
      a_{11,\pm} & a_{12,\pm}\\
      a_{12,\pm} & a_{22,\pm}
    \end{array}
  \right]
\end{equation}
be their inverse square roots,
and 
\begin{equation}
  \label{eq:qpm}
  Q_{\pm} = \frac{1}{\alpha_{\pm}}\left[
    \begin{array}{cc}
      a_{11,\pm} & a_{12,\pm}\\
      -a_{12,\pm} & a_{11,\pm}\\
    \end{array}
    \right]
\end{equation}
be two related orthogonal matrices, where
$\alpha_{\pm}=\sqrt{a_{11,\pm}^2+a_{12,\pm}^2}>0$. To uniquely determine $u^{\rm
  tot}$, we need a proper radiation condition at infinity, and this relies on a
precise definition of outgoing waves, as shown below.
\begin{mydef}
  \label{def:og}
  Let $B_R=\{x:|x|<R\}$ for any $R>0$. A wave field $u(x)$ is said to be
  outgoing in $\Omega^{\pm}$ if for sufficiently large $R$, $u\in
  C^{2}(\Omega^{\pm}\backslash \overline{B_R})$ satisfies (\ref{eq:Hhel}) on
  $\Omega^{\pm}\backslash \overline{B_R}$ and (\ref{eq:cont:cond}) across
  $\Gamma\backslash \overline{B_R}$, and if $U_\pm(X) = u(M_{\pm}^{1/2}Q_{\pm}^{T}X)$
  satisfy the following half-plane SRC (hpSRC)
\begin{align}
  \label{eq:src}
  \lim_{|X|\to\infty}\sqrt{|X|}\left( \frac{\partial}{\partial |X|}-\bi k_0\right)U_{\pm}(X)=0,\quad \pm X_2>0,
\end{align}
{uniformly in all directions $X/|X|\in\{(\cos\beta,\sin\beta):0\leq \pm\beta\leq \pi\}$}, where 
$X=(X_1,X_2)^{T}$ indicates a new coordinate system via $X=Q_{\pm}M_{\pm}^{-1/2}x$.
\end{mydef}
The outgoing behavior of $u$, if satisfying (\ref{eq:src}), shall be justified
later in Section 4. We make some remarks below.
\begin{myremark}
  \label{rem:src}
  The transition matrices $Q_{\pm}M_{\pm}^{-1/2}$ can directly transform the
  anisotropic Helmholtz equation (\ref{eq:Hhel}) on $\Omega^{\pm}\backslash \overline{B_R}$
  to the following isotropic Helmholtz equation
  \begin{equation}
    \label{eq:helm:diffcord}
    \Delta_{X}U_\pm(X) + k_0^2U_\pm(X) = 0,
  \end{equation}
  where $\Delta_{X}=\sum_{j=1}^2\partial^2_{X_j}$. Therefore, it is natural
  to use the hpSRC (\ref{eq:src}) in the $X$-coordinate system to characterize
  outgoing waves in $\Omega^{\pm}$. Note that the two matrices
  $Q_{\pm}M_{\pm}^{-1/2}$ are different so that the hpSRC
  (\ref{eq:src}) holds in two different coordinate systems.
\end{myremark}
\begin{myremark}
  \label{rem:Q}
  The orthogonal matrices $Q_{\pm}$ in fact can be replaced by any two other
  orthogonal matrices to keep the resulting Helmholtz equations isotropic.
  Nevertheless, we choose (\ref{eq:qpm}) since $Q_{\pm}M_{\pm}^{-1/2}$ are upper
  triangular so that the horizontal $x_2$-axis is exactly the horizontal
  $X_2$-axis. Thus, $c\cdot Q_{\pm}M_{\pm}^{-1/2}$ for any $c>0$ can be an
  acceptable choice. As we shall see, such a choice of $Q_{\pm}$ is essential in
  computing the background Green's functions and in the setup of an RPML.
\end{myremark}
Certainly, it is inaccurate to directly assume that $u^{\rm tot}$ is outgoing
without considering the incident part of $u^{\rm tot}$. For simplicity, we shall
assume that the incident wave is specified in $\Omega^+$ only and that $M_+$ is
the $2\times 2$ identity matrix $I_2$ so that $Q_+=I_2$. According to Remark~\ref{rem:src}, if
$M_+\neq I_2$, we use $Q_+M_+^{-1/2}x$ as the new two-dimensional (2D)
coordinate system in the $x_1Ox_2$-plane so that $u^{\rm tot}$ satisfies
(\ref{eq:Hhel}) but with $M(x)$ replaced by $Q_+M_+^{-1/2}M(x) M_+^{-1/2}Q_+^T$.
Thus, we can specify incident waves of simple forms in $\Omega^+$. From now on,
we shall always identify $X$ as $Q_-M_-^{-1/2}x$ and $x$ as $M_-^{1/2}Q_-^{T}X$
when $x\in\Omega^-$, and for a generic domain $\Omega\subset\Omega^-$, we shall
call $\Omega_N$ the image of $\Omega$ if
$\Omega_N=\{X=Q_-M_-^{-1/2}x:x\in\Omega\}$, and hence shall call $\Omega$ the
preimage of $\Omega_N$.

In this paper, we consider only plane and cylindrical incident waves, and they
are separately discussed below. For a plane incident wave $u^{\rm
  inc}(x;\theta)=e^{\bi k_0 (\cos\theta x_1 - \sin\theta x_2)},
\theta\in(0,\pi),$ specified in $\Omega^+$, we distinguish two cases. If
$\Gamma$ is the unperturbed straight line $x_2=0$, then the total field,
referred to as the background solution $u^{\rm tot}_{\rm b}(x;\theta)$, and can
be predetermined by the method of Fourier transform; its closed form is
presented in (\ref{eq:utotb:pw1}) and (\ref{eq:utotb:pw2}) in Section 3.2. In
general, if $\Gamma$ is a locally-perturbed straight line, we enforce the
following radiation condition
\begin{center}
  (RC1):\quad\quad $u^{\rm tot}-u^{\rm tot}_{\rm b}$ is outgoing in $\Omega^{\pm}$.
\end{center}
For a cylindrical incident wave $u^{\rm inc}(x;x^*)=\frac{\bi}{4}
H_0^{(1)}(k_0|x-x^*|)$ where $x^*=(x_1^*,x_2^*)^{T}\in\Omega^+$ denotes the
exciting source point, the right-hand side of (\ref{eq:Hhel}) should be replaced
by $-\delta(x-x^*)$. Thus, $u^{\rm tot}(x;x^*)$ represents the Green's function
excited by the source point $x^*$. The background solution, the Green's function
$u^{\rm tot}_b(x;x^*)$, for the unperturbed case $\Gamma=\{x:x_2=0\}$, can again
be predetermined; see Section 3.1 for details. In the locally perturbed case,
one can still enforce (RC1), but it is practically more efficient to enforce the
following radiation condition
\begin{center}
(RC2):\quad\quad$u^{\rm tot}$ is outgoing in $\Omega^{\pm}$,
\end{center}
since $u_b^{\rm tot}$ is not involved. The equivalence of (RC1) and (RC2) for
cylindrical-wave incidences shall be justified in Corollary~\ref{cor:equivRC12}.

The objective of this paper is to develop an efficient numerical method to
compute $u^{\rm tot}$ governed by (\ref{eq:Hhel}), (\ref{eq:cont:cond}), and one
of the two radiation conditions (RC1) and (RC2) depending on the type of the
incidence. To achieve this, it is clear that the background solutions $u^{\rm
  tot}_{\rm b}$ should be computed in advance, and this is the main content of
the next section.
\section{Background solutions}
In this section, we assume $\Gamma=\{x\in\mathbb{R}^2:x_2=0\}$ so that
$\Omega_\pm=\mathbb{R}_\pm^2:=\{x\in\mathbb{R}^2:\pm x_2>0\}.$ We shall use the
method of Fourier transform to compute the background solutions $u^{\rm tot}_b$
for plane and cylindrical incident waves.
\subsection{Cylindrical incident wave}
By convention, it is more appropriate to use $G(x;x^*)$ instead of $u^{\rm
  tot}_b(x;x^*)$ to denote the background Green's function. We consider first
the case $x^*\in\mathbb{R}_+^2$, i.e., $x_2^*>0$. Recall that
$\alpha_-=\sqrt{a_{11,-}^2+a_{12,-}^2}$, and that we have assumed $M_+=I_2$.
According to Remark \ref{rem:src}, $G(x;x^*)$ satisfies
\begin{align}
\label{eq:green1}
  \Delta G\left(x;x^{*}\right)+k_0^{2}G\left(x;x^{*}\right)=&-\delta\left(x-x^{*}\right),\quad x_2>0,\\
\label{eq:green2}
  \Delta_{X} G_-\left(X;x^{*}\right)+k_0^{2}G_{-}\left(X;x^{*}\right)=&0,\quad\quad\quad\quad\quad\quad X_2<0,\\
\label{eq:green3}
  G(x;x^*) = G_-(X;x^*),\quad \partial_{x_2} G(x;x^*) =& \alpha_-\sqrt{|M_-|}\partial_{X_2}G_-(X;x^*),\quad x_2=0=X_2,
\end{align}
where $G_-(X;x^*)=G(M_-^{1/2}Q^T_-X;x^*)$, $|M_-|$ denotes the determinant of
$M_-$, (\ref{eq:green3}) is derived from (\ref{eq:cont:cond}). Note that $x_2<0$
and $X_2<0$ represent the same lower-half plane $\mathbb{R}_-^2$, i.e., the
image of $\mathbb{R}_-^2$ is itself. Let
\begin{align*}
  \hat{G}\left(x_{2}; x^*,\xi\right)=\int_{-\infty}^{\infty} G(x;x^{*}) e^{\mathbf{i} \xi x_{1}} \mathrm{d} x_{1} \quad{\rm and}\quad \hat{G}_{-}\left(X_{2}; x^*,\xi\right)=\int_{-\infty}^{\infty} G_{-}(X;x^{*}) e^{\mathbf{i} \xi X_{1}} \mathrm{d} X_{1}
\end{align*}
be the one-dimensional Fourier transforms of $G$ and $G_-$ w.r.t $x_1$ and $X_1$
variables, respectively. The governing ordinary differential equations of
$\hat{G}$ and $\hat{G}_-$ are
\begin{align}
\label{eq:gov:hG}
\hat{G}\dd\left(x_{2} ; x^*,\xi\right)+\mu(\xi)^2\hat{G}\left(x_{2} ; \xi\right)=&-e^{\mathbf{i}\xi x_{1}^{*}}\delta\left(x_2-x_2^{*}\right),\quad x_2>0,\\ 
\label{eq:gov:hG-}
\hat{G}_{-}^{\prime\prime}\left(X_{2} ;x^*, \xi\right)+\mu(\xi)^2\hat{G}_{-}\left(X_{2} ;x^*, \xi\right)=&0,\quad X_2<0,
\end{align}
where $\mu(\xi)= \sqrt{k_0^2-\xi^2}$. Throughout this paper, we use the negative
real axis as the branch cut of the square-root function $\sqrt{\cdot}$ to limit
its argument onto $(-\pi/2,\pi/2]$. Thus, we seek $\hG$ and $\hG_-$ in the form
of
\begin{align}
  \hG\left(x_{2} ; x^*,\xi\right) =& \frac{\bi}{2\mu(\xi)}e^{\bi \mu(\xi) |x_2-x_2^*|+\bi \xi x_1^*} + A(\xi) e^{\bi \mu(\xi) x_2},\\
  \hG_-(X_2;x^*,\xi) =& B(\xi) e^{-\bi \mu(\xi) X_2},
\end{align}
where we have disregarded the downgoing wave $e^{-\bi\mu x_2}$ in $\hG$ and the
upgoing wave $e^{\bi \mu X_2}$ in $\hG_-$. The continuous condition
(\ref{eq:green3}) implies
\begin{equation}
\label{eq:hcond}
  \hG_-(0;x^*,\xi)=\alpha_-\hG(0;x^*,\alpha_-\xi),\quad\hG_-'(0;x^*,\xi)=|M_-|^{-1/2}\hG'(0;x^*,\alpha_-\xi),
\end{equation}
so that
\begin{align*}
  A(\alpha_-\xi) + \frac{\bi}{2\mu(\alpha_-\xi)}e^{\bi \mu(\alpha_-\xi) x_2^*+\bi\alpha_-\xi x_1^*} =& \frac{1}{\alpha_-}B(\xi),\\
  A(\alpha_-\xi) -\frac{\bi}{2\mu(\alpha_-\xi)}e^{\bi \mu(\alpha_-\xi) x_2^*+\bi\alpha_-\xi x_1^*} =& -\frac{|M_-|^{1/2}\mu(\xi)}{\mu(\alpha_-\xi)}B(\xi).
\end{align*}
Solving the above linear system gives rise to
\begin{align}
  A(\xi) =& \frac{\bi(\mu(\xi)-|M_-|^{1/2}\alpha_-\mu(\alpha_-^{-1}\xi))}{2\mu(\xi)(\mu(\xi)+|M_-|^{1/2}\alpha_-\mu(\alpha_-^{-1}\xi))}e^{\bi \mu(\xi) x_2^*+\bi\xi x_1^*},\\
  B(\xi) =&\frac{\alpha_-\bi}{\mu(\alpha_-\xi)+|M_-|^{1/2}\alpha_-\mu(\xi)}e^{\bi \mu(\alpha_-\xi) x_2^*+\bi\alpha_-\xi x_1^*}.
\end{align}
Thus, inverse Fourier transforming $\hat{G}$, we obtain for
$x\in\mathbb{R}^2_+$,
\begin{align}
  \label{eq:G:++}
  G(x;x^*) 
  =&\Phi(x;x^*) + \frac{\bi}{4\pi}\int_{-\infty}^{\infty}\frac{\mu(\xi)-|M_-|^{1/2}\alpha_-\mu(\alpha_-^{-1}\xi)}{\mu(\xi)(\mu(\xi)+|M_-|^{1/2}\alpha_-\mu(\alpha_-^{-1}\xi))}e^{\bi \mu(\xi) (x_2^*+x_2)+\bi\xi (x_1^*-x_1)}d\xi,
\end{align}
where $\Phi(x;x^*)=\frac{\bi}{4} H_0^{(1)}(k_0|x-x^*|)$ appears due to the
identity
\[
H_0^{(1)}(k_0|x-x^*|) = \frac{1}{\pi}\int_{-\infty}^{\infty}\frac{1}{\mu(\xi)}e^{\bi \mu(\xi) |x_2-x_2^*|+\bi \xi (x_1^*-x_1)}d\xi.
\]
Similarly, for $x\in\mathbb{R}_-^2$, 
\begin{align}
  \label{eq:G:+-}
  G(x;x^*)=G_-(X;x^*) = \frac{\alpha_-\bi}{2\pi}\int_{-\infty}^{\infty}\frac{e^{\bi \mu(\alpha_-\xi) x_2^*-\bi \mu(\xi)X_2+\bi\alpha_-\xi x_1^*-\bi\xi X_1}}{\mu(\alpha_-\xi)+|M_-|^{1/2}\alpha_-\mu(\xi)}d\xi.
\end{align}
Note that the denominators
\[
  \mu(\alpha_-\xi)+|M_-|^{1/2}\alpha_-\mu(\xi)\quad{\rm and}\quad\mu(\xi)+|M_-|^{1/2}\alpha_-\mu(\alpha_-^{-1}\xi)
\]
vanish only when $\alpha_-=1$ and $\xi=\pm k_0$, in which case the above
integrals still exist as Riemann integrals. 

For completeness, we give the closed form of $G$ for $x^*\in \mathbb{R}_-^2$.
By similar derivations, for $x\in\mathbb{R}^2_+$, 
\begin{align}
  \label{eq:G:-+}
  G(x;x^*)&=\frac{\bi}{2\pi}\int_{-\infty}^{\infty}\frac{e^{-\bi \mu(\alpha^{-1}_-\xi) X_2^*+\bi\alpha_-^{-1}\xi X_1^*+\bi\mu(\xi)x_2-\bi\xi x_1}}{|M_-|^{1/2}\alpha_-\mu(\alpha_-^{-1}\xi)+\mu(\xi)}d\xi,
\end{align}
and for $x\in\mathbb{R}_-^2$,
\begin{align}
  \label{eq:G:--}
  G(x;x^*)=&G_-(X;x^*)=|M_-|^{-1/2}\Phi(X;X^*) \nonumber\\
  &+ \frac{\bi}{4\pi|M_-|^{1/2}}\int_{-\infty}^{\infty}\frac{\alpha_-\mu(\xi)|M_-|^{1/2}-\mu(\alpha_-\xi)}{\mu(\xi)(\alpha_-\mu(\xi)|M_-|^{1/2}+\mu(\alpha_-\xi))}e^{-\bi \mu(\xi) (X_2^*+X_2)+\bi\xi (X_1^*-X_1)}d\xi,
\end{align}
where $X^*=Q_-M_-^{-1/2}x^*$. From the above formulae, we see that $G$ satisfies
the following reciprocity relation
\begin{equation}
  \label{eq:rec:rel}
  G(x;x^*) = G(x^*;x).
\end{equation}
The following lemma describes the asymptotic behavior of $G(x;x^*)$ as
$|x|\to\infty$.
\begin{mylemma}
  \label{lem:asym:G}
  For any $x^*\in\mathbb{R}^2$ with $x_2^*\neq 0$, the background Green's
  function $G(x;x^*)$
  has the following properties:\\
  (a). Let $x^*\in\mathbb{R}_+^2$. For $x_2\geq 0$,
      \begin{align}
        \label{eq:fp:++}
        G(x;x^*) = &\frac{e^{\bi k_0|x|}}{\sqrt{|x|}}\left[G_{\infty}^{++}(\beta;x^*) + {\cal O}(|x|^{-1})\right], 
      \end{align}
      as $|x|\to\infty$, where $\beta\in[0,\pi]$ is such that
      $x=(|x|\cos\beta,|x|\sin\beta)$ and
      \begin{align*}
        G_{\infty}^{++}(\beta;x^*)=&\frac{e^{\bi \pi/4-\bi k_0\sin\beta x_2^*-\bi k_0\cos\beta x_1^*}}{\sqrt{8\pi k_0}}\nonumber\\
        &+\frac{e^{\bi \pi/4+\bi k_0\sin\beta x_2^*-\bi k_0\cos\beta
            x_1^*}}{\sqrt{8\pi k_0}}
        \frac{k_0\sin\beta-|M_-|^{1/2}\alpha_-\mu(-\alpha_-^{-1}k_0\cos\beta)}{k_0\sin\beta+|M_-|^{1/2}\alpha_-\mu(-\alpha_-^{-1}k_0\cos\beta)}.
      \end{align*}
      For $x_2\leq 0$,
      \begin{align}
        \label{eq:fp:-+}
        G_-(X;x^*) = \frac{e^{\bi k_0|X|}}{\sqrt{|X|}}\left[G_{\infty}^{-+}(\tilde{\beta};x^*)+ {\cal O}(|X|^{-1})\right], 
      \end{align}
      as $|X|\to\infty$, where $\tilde{\beta}\in[\pi,2\pi]$ is such that
      $X=(|X|\cos\tilde{\beta},|X|\sin\tilde{\beta})$ and 
      \[
        G_{\infty}^{-+}(\tilde{\beta};x^*)=\frac{k_0\sin\tilde{\beta}\alpha_-e^{\bi
            \pi/4+\bi \mu(\alpha_-k_0\sin\tilde{\beta}) x_2^*-\bi
            \alpha_-k_0\cos\tilde{\beta} x_1^*}}{\sqrt{2\pi
            k_0}(|M_-|^{1/2}\alpha_-k_0\sin\tilde{\beta}+\mu(\alpha_-k_0\cos\tilde{\beta}))}.
      \]\\
      (b). Let $x^*\in\mathbb{R}_-^2$. For $x_2\geq 0$,
      \begin{align}
        \label{eq:fp:+-}
        G(x;x^*) = \frac{e^{\bi k_0|x|}}{\sqrt{|x|}}\left[G_{\infty}^{+-}(\beta;x^*)
        + {\cal O}(|x|^{-1})\right], 
      \end{align}
      as $|x|\to\infty$, where
      \[
G_{\infty}^{+-}(\beta;x^*)=\frac{k_0\sin{\beta}e^{\bi \pi/4-\bi \mu(\alpha_-k_0\sin{\beta})  X_2^*-\bi \alpha_-^{-1}k_0\cos{\beta} X_1^*}}{\sqrt{2\pi k_0}(k_0\sin\beta + |M_-|^{1/2}\alpha_-\mu(\alpha_-^{-1}k_0\cos{\beta}))},
      \]
      and $X^*=(X_1^*,X_2^*)^{T}$. For $x_2\leq 0$,
      \begin{align}
        \label{eq:fp:--}
        G_-(X;x^*) = &\frac{e^{\bi k_0|X|}}{\sqrt{|X|}}\left[G_{\infty}^{--}(\tilde{\beta};x^*) + {\cal O}(|X|^{-1})\right], 
      \end{align}
      as $|X|\to\infty$, where
      \begin{align*}
G_{\infty}^{--}(\tilde{\beta};x^*)=&\frac{e^{\bi \pi/4-\bi k_0\sin\tilde{\beta} X_2^*-\bi k_0\cos\tilde{\beta} X_1^*}}{\sqrt{8|M_-|\pi k_0}}\nonumber\\ 
                     &+\frac{e^{\bi \pi/4-\bi k_0\sin\tilde{\beta} X_2^*-\bi k_0\cos\tilde{\beta} X_1^*}}{\sqrt{8|M_-|\pi k_0}} \frac{|M_-|^{1/2}\alpha_-k_0\sin\tilde{\beta}-\mu(\alpha_-k_0\cos\tilde{\beta})}{|M_-|^{1/2}\alpha_-k_0\sin\tilde{\beta}+\mu(\alpha_-k_0\cos\tilde{\beta})}.
      \end{align*}
      \\
      (c). $G(x;x^*)$ satisfies (RC2).\\
      In the above, the prefactors in the ${\cal O}$-terms do not depend on
      $\beta$ or $\tilde{\beta}$.
  \begin{proof}
    Based on contour deformations similar to those in the proof of Lemma~2.1 in
    \cite{lu21}, it is straightforward to verify the above properties.
  \end{proof}
\end{mylemma}
By convention, $G_{\infty}^{\pm,\pm}$ in the brackets of
(\ref{eq:fp:++})-(\ref{eq:fp:--}) constitute the far-field patterns of
$G(x;x^*)$. Lemma~\ref{lem:asym:G}(c) implies
\begin{mycorollary}
  \label{cor:equivRC12}
  For any cylindrical incident wave $u^{\rm inc}(x;x^*)$ with
  $x^*\in\Omega^{\pm}$ and any locally perturbed straight line $\Gamma$, (RC1)
  and (RC2) for the total field $u^{\rm tot}(x;x^*)$ are equivalent.
\end{mycorollary}
\subsection{Plane incident wave}
Suppose now $u^{\rm inc}(x;\theta)=e^{\bi k_0 (\cos\theta x_1 - \sin\theta
  x_2)}$ for $\theta\in(0,\pi)$. Then, the background solution $u^{\rm
  tot}_b(x;\theta)$ satisfies
\begin{align}
\label{eq:utotb1}
  \Delta u^{\rm tot}_b\left(x;\theta\right)+k_0^{2}u^{\rm tot}_b\left(x;\theta\right)=&0,\quad x_2>0,\\
\label{eq:utotb2}
  \Delta_{X} U^{\rm tot}_{b,-}\left(X;\theta\right)+k_0^{2}U^{\rm tot}_{b,-}\left(X;\theta\right)=&0,\quad X_2<0,\\
\label{eq:utotb3}
  u^{\rm tot}_b(x;\theta) = U^{\rm tot}_{b,-}(X;\theta),\quad \partial_{x_2} u^{\rm tot}_b(x;\theta) =& \alpha_-\sqrt{|M_-|}\partial_{X_2}U^{\rm tot}_{b,-}(X;\theta),\quad x_2=0=X_2,
\end{align}
where $U^{\rm tot}_{b,-}(X;\theta) = u^{\rm tot}_b(Q_-M_-^{-1/2}X;\theta)$ for
$x_2<0$. By the same approach as above, we obtain
\begin{align}
  \label{eq:utotb:pw1}
u^{\rm tot}_b(x;\theta)=& e^{\bi k_0 (\cos\theta x_1 -
                          \sin\theta x_2)} \nonumber\\
  &+ \frac{k_0\sin\theta - |M_-|^{1/2}\mu(\alpha_-^{-1}k_0\cos\theta)\alpha_-}{k_0\sin\theta+|M_-|^{1/2}\mu(\alpha_-^{-1}k_0\cos\theta)\alpha_-}e^{\bi k_0(\cos\theta x_1+\sin\theta x_2)},\quad x\in\mathbb{R}_+^2,\\
  \label{eq:utotb:pw2}
U^{\rm tot}_{b,-}(X;\theta)=&\frac{2k_0\sin\theta e^{-\bi \mu(\alpha_-^{-1}k_0\cos\theta)X_2+\bi \alpha_{-}^{-1}k_0\cos\theta X_1} }{k_0\sin\theta + |M_-|^{1/2}\alpha_-\mu(\alpha_-^{-1}k_0\cos\theta)},\quad X\in\mathbb{R}_-^2.
\end{align}
Note that to obtain the above formulae, we have implicitly assumed that the
reflected and transmitted waves are upgoing and downgoing, respectively.

\section{Green's representation formula}
In this section, we shall use the background Green's function to derive Green's
representation formula for a locally perturbed straight line $\Gamma$. Let
$\Gamma_T$ denote the boundary of a bounded Lipschitz domain $\Omega_T$
enclosing the perturbed part of $\Gamma$, and let $\Gamma_T^{\pm} =
\Gamma_T\cap\Omega^{\pm}$. Define
\[
  \nu_c(x) = \begin{cases}
    \nu(x), & x\in\Gamma_T^+,\\
    M_-\nu(x), & x\in\Gamma_T^-,
    \end{cases}
\]
as the conormal vector along $\Gamma_T$, where $\nu(x)$ is the outer unit
normal vector of $\Gamma_T$. Let $\Gamma_{N,T}^{-}$ be the image
of the lower part $\Gamma_T^-$. We have the following theorem regarding Green's
representation formula.
\begin{mytheorem}
  \label{thm:greenrep}
  Let $B_r\subset \Omega_T$ for some sufficiently large $r>0$. A function $u$ is
  outgoing in $\Omega^{\pm}$ if and only if for
  $x\in\mathbb{R}^2\backslash\overline{\Omega_T}$,
\begin{align}
  \label{eq:rep:G}
u(x) =&\int_{\Gamma_T}\left[\partial_{\nu_c(y)}G(y;x) u(y)  -G(y;x)\partial_{\nu_c}u(y) \right]ds(y),
\end{align}
or alternatively, 
\begin{align}
  \label{eq:rep:GG-}
u(x) =&\int_{\Gamma_T^+}\left[\partial_{\nu(y)}G(y;x) u(y) -G(y;x)\partial_{\nu}u(y)  \right]ds(y)\nonumber\\
&+|M_-|^{1/2}\int_{\Gamma_{N,T}^{-}}\left[\partial_{\nu^N(Y)}G_-(Y;x) U_-(Y) - G_-(Y;x)\partial_{\nu^N}U_-(Y)\right]ds(Y),
\end{align}
where we recall $U_-(X)=u(Q_-M_-^{-1/2}X)$ and $G_-(Y;x)=G(Q_-M_-^{1/2}Y;x)$,
$\nu$ denotes the outer unit normal vector of $\Gamma_T^{+}$, and $\nu^N$
denotes the outer unit normal vector of $\Gamma_{N,T}^{-}$.
\begin{proof}
  We consider the ``only if'' part as the other part is straightforward by
  Lemma~\ref{lem:asym:G}. Without loss of generality, we assume $x\in\Omega^+$.
  For sufficiently large $R>0$, let $\Gamma_R^+=\partial B_R\cap \Omega^+$ and
  let $\Gamma_R^-$ be the preimage of the lower-half circle
  $\Gamma_{N,R}^{-}=\partial B_{\alpha_-R}\cap \mathbb{R}_-^2$. We choose the
  boundary $\partial B_{\alpha_-R}$ since $\Gamma_R=\Gamma_R^+\cup
  \overline{\Gamma_R^-}$ forms a closed curve enclosing $\Gamma_T$. Let
  $\Gamma_{TR}$ be the union of the two line segments on $\Gamma$ between
  $\Gamma_T$ and $\Gamma_R$, and $\Gamma_{N,TR}$ be its image. Green's third
  identity implies
  \begin{align*}
    u(x) = \int_{\Gamma_R^+\cup\Gamma_T^+\cup \Gamma_{TR}}\left[G(y;x)\partial_{\nu(y)}u(y) - \partial_{\nu(y)}G(y;x) u(y)  \right]ds(y).
  \end{align*}
  On $\Gamma_{TR}$,
  $ds(y)=dy_1=\alpha_-^{-1}dY_1=\alpha_-^{-1}ds(Y)$ and
  \[
    \partial_{\nu(y)}=-\partial_{y_2}=-\alpha_-|M_-|^{1/2}\partial_{Y_2}=-\alpha_-|M_-|^{1/2}\partial_{\nu^N(Y)},
  \]
  so that
  \begin{align*}
    &\int_{\Gamma_{TR}}\left[G(y;x)\partial_{\nu(y)}u(y) - \partial_{\nu(y)}G(y;x) u(y)  \right]ds(y)\\
    =&-|M_-|^{1/2}\int_{\Gamma_{N,TR}}\left[G_-(Y;x)\partial_{\nu^N(Y)}U_-(Y) - \partial_{\nu^N(Y)}G(Y;x) U_-(Y)  \right]ds(Y)\\
    =&|M_-|^{1/2}\int_{\Gamma_{N,T}^{-}\cup\Gamma_{N,R}^{-}} \left[G_-(Y;x)\partial_{\nu^N(Y)}U_-(Y) - \partial_{\nu^N(Y)}G(Y;x) U_-(Y)  \right]ds(Y).
  \end{align*}
  Taking advantage of the hpSRC (\ref{eq:src}) in Definition~\ref{def:og} and
  Lemma~\ref{lem:asym:G}, we obtain
  \begin{align*}
    &\int_{\Gamma_R^+}\left[G(y;x)\partial_{\nu(y)}u(y) - \partial_{\nu(y)}G(y;x) u(y)  \right]ds(y) \\
    +& |M_-|^{1/2}\int_{\Gamma_{N,R}^{-}} \left[G_-(Y;x)\partial_{\nu^N(Y)}U_-(Y) - \partial_{\nu^N(Y)}G(Y;x) U_-(Y)  \right]ds(Y)\to 0,\quad{\rm as}\quad R\to\infty,
  \end{align*}
  by arguing as in the proof of Theorem 2.5 in \cite{colkre13}. Consequently,
  (\ref{eq:rep:GG-}) holds. As for (\ref{eq:rep:G}), it is straightforward to
  verify that $ds(y)\partial_{\nu_c(y)}=|M_-|^{1/2}ds(Y)\partial_{\nu^N(Y)}$ under 
  the change of variable $Y=Q_-M_-^{-1/2}y$.
\end{proof}
\end{mytheorem}
Clearly, Theorem~\ref{thm:greenrep} and Lemma~\ref{lem:asym:G} explain the
reasonability of Definition~\ref{def:og} regarding outgoing waves since $G$ is
outgoing in $\Omega^{\pm}$. To conclude this section, we discuss two byproducts
of Green's representation formulae (\ref{eq:rep:G}) and (\ref{eq:rep:GG-}).
Firstly, the far-field pattern of any outgoing wave can now be well defined.
According to the far-field pattern of $G$ shown in Lemma~\ref{lem:asym:G}, we
see immediately from (\ref{eq:rep:GG-}) that any outgoing wave field $u$
asymptotically behaves as follows: for
$x=|x|(\cos\beta,\sin\beta)\in\mathbb{R}_+^2$ with $\beta\in[0,\pi]$,
\begin{align}
  \label{eq:fp:u+}
  u(x) = &\frac{e^{\bi k_0|x|}}{\sqrt{|x|}}\left[\hat{u}^+_{\infty}(\beta)+{\cal O}\left(\frac{1}{|x|}\right)\right], 
\end{align}
as $|x|\to\infty$, where 
\begin{align*}
\hat{u}_{\infty}^+(\beta)=&\int_{\Gamma_T^+}\left[\partial_{\nu(y)}G_{\infty}^{++}(\beta;y) u(y) -G_{\infty}^{++}(\beta;y)\partial_{\nu}u(y)  \right]ds(y)\nonumber\\
                          &+|M_-|^{1/2}\int_{\Gamma_{N,T}^{-}}\left[\partial_{\nu^N(Y)}G_{\infty}^{+-}(\beta;Y) U_-(Y) - G_{\infty}^{+-}(\beta;Y)\partial_{\nu^N}U_-(Y)\right]ds(Y);
\end{align*}
for $x\in\mathbb{R}_-^2$ so that
$X=Q_-M_-^{-1/2}x=|X|(\cos\tilde{\beta},\sin\tilde{\beta})$ with
$\tilde{\beta}\in[\pi,2\pi]$,
\begin{align}
  \label{eq:fp:u-}
  u(x) = &\frac{e^{\bi k_0|X|}}{\sqrt{|X|}}\left[\hat{u}^-_{\infty}(\tilde{\beta})+{\cal O}\left(\frac{1}{|X|}\right)\right], 
\end{align}
as $|X|\to\infty$, where 
\begin{align*}
\hat{u}_{\infty}^-(\tilde{\beta})=&\int_{\Gamma_T^+}\left[\partial_{\nu(y)}G_{\infty}^{-+}(\tilde{\beta};y) u(y) -G_{\infty}^{-+}(\tbeta;y)\partial_{\nu}u(y)  \right]ds(y)\nonumber\\
                          &+|M_-|^{1/2}\int_{\Gamma_{N,T}^{-}}\left[\partial_{\nu^N(Y)}G_{\infty}^{--}(\tbeta;Y) U_-(Y) - G_{\infty}^{--}(\tbeta;Y)\partial_{\nu^N}U_-(Y)\right]ds(Y).
\end{align*}
Consequently, $\hat{u}^{+}_{\infty}(\beta)$ and
$\hat{u}^-_{\infty}(\tilde{\beta})$ can be defined as the far-field pattern of
the outgoing field $u$.

Secondly, the unbounded domain can now be truncated onto $\Omega_T$. Let $x$ approach $\Gamma_T$ in (\ref{eq:rep:G}), we obtain the following
transparent boundary condition (TBC)
\begin{equation}
  \label{eq:tbc}
  ({\cal I} - {\cal K}_b) [u](x)= -{\cal S}_b[\partial_{\nu_c}u](x)
\end{equation}
on $\Gamma_T$, where ${\cal I}$ is the identity operator such that ${\cal
  I}[u]=u$, and ${\cal K}_b$ and ${\cal S}_b$ are two boundary integral operators on $\Gamma_T$ defined as
follows: for any $x\in\Gamma_T$ and $\psi\in C^{\infty}(\Gamma_T)$,
\begin{align}
  \label{eq:int:S}
  {\cal S}_b[\psi](x):=&2\int_{\Gamma_T} G(y;x)\psi(y)ds(y),\\
  \label{eq:int:K}
  {\cal K}_b[\psi](x):=&2{\rm p.v.}\int_{\Gamma_T} \partial_{\nu_c(y)}G(y;x)\psi(y)ds(y),
\end{align}
where ${\rm p.v.}$ indicates Cauchy principal value. In our scattering problem,
according to the radiation conditions (RC1) and (RC2), $u$ can represent $u^{\rm
  tot}-u^{\rm tot}_b$ or $u^{\rm tot}$ depending on the type of the incident
wave.

Theoretically, the TBC (\ref{eq:tbc}) plays a central role in proving the
well-posedness of the scattering problem (\ref{eq:Hhel}) and
(\ref{eq:cont:cond}) under the radiation condition (RC1) or (RC2), but
unfortunately, standard approaches coupling TBCs to variational formulations
\cite{hulurat21} break down here. This is because that the integral operators
${\cal S}_b$ and ${\cal K}_b$ lose properties such as strong ellipticity and
compactness. In a subsequent work \cite{lu21ort}, we shall present a new
framework to establish the well-posedness theory.

Numerically, one may incorporate any standard numerical method with
(\ref{eq:tbc}) to compute $u^{\rm tot}$ on the bounded domain $\Omega_T$ by
properly discretizing the two integral operators ${\cal S}_b$ and ${\cal K}_b$.
Nevertheless, it becomes essential to develop a fast and accurate algorithm to
evaluate $G$ and its derivatives. If one does not wish to use the TBC
(\ref{eq:tbc}) in a numerical method, then an artificial boundary condition must
be developed, as shall be discussed in the next section.

\section{Regionalized PML and BIE method}
In this section, we shall propose a stable RPML to truncate $\Omega^{\pm}$, and
shall develop a high-accuracy BIE method to
numerically compute $u^{\rm tot}$. Unless otherwise indicated, we shall assume
the incident wave to be the plane wave $u^{\rm inc}(x;\theta)$. For simplicity,
we shall suppress the argument $\theta$.
\subsection{Governing equations}
Let $u^{\rm og}_{\pm} = u^{\rm og}|_{\Omega^{\pm}}$ be the two outgoing waves in
$\Omega^{\pm}$ according to (RC1). Then, $u^{\rm og}_{\pm}$ satisfy
\begin{align}
  \label{eq:uog1}
  \nabla\cdot(M_{\pm} \nabla u_{\pm}^{\rm og}) + k_0^2 u_{\pm}^{\rm og} =& 0,\quad{\rm on}\quad\Omega^{\pm},\\
  \label{eq:uog2}
  u_{+}^{\rm og}(x) - u_{-}^{\rm og}(x) =& \lim_{y\to x^-}u^{\rm tot}_b(y)-\lim_{y\to x^+}u^{\rm tot}_b(y),\quad x\in \Gamma,\\
  \label{eq:uog3}
  \partial_{\nu_c}u_{+}^{\rm og}(x) - \partial_{\nu_c}u_{-}^{\rm og}(x) =& \lim_{y\to x^-}\partial_{\nu_c}u^{\rm tot}_b(y)-\lim_{y\to x^+}\partial_{\nu_c}u^{\rm tot}_b(y),\quad x\in \Gamma,
\end{align}
where we now set $\nu(x)$ to be the unit normal vector of $\Gamma$ at $x$
pointing towards $\Omega^-$, $\nu_c=M\nu$, and $x^\pm$ indicate one-sided limits
taken from $\Omega^\pm$, respectively. Note that $u^{\rm tot}_b(x)$ should be
defined as (\ref{eq:utotb:pw1}) for $x\in\Omega^+$ and (\ref{eq:utotb:pw2}) for
$x\in\Omega^-$ even if $x$ is outside one of the two domains $\mathbb{R}_\pm^2$,
so that (\ref{eq:uog1}) holds.

Recall the assumption $M_+=I_2$ and the change of coordinates
$X=Q_-M_-^{-1/2}x$. Let $U_-^{\rm og}(X)=u_-^{\rm og}(M_-^{1/2}Q_-^{T}X)$,
$\Omega^-_N$ be the image of $\Omega^-$, and $\Gamma_N=\partial\Omega_{N}^-$.
Then, (\ref{eq:uog1})-(\ref{eq:uog3}) become
\begin{align}
  \label{eq:uog1:n}
  \Delta u^{\rm og}_+(x) + k_0^2 u^{\rm og}_+(x) =& 0,\quad x\in\Omega^+,\\
  \label{eq:uog2:n}
  \Delta_X U^{\rm og}_- + k_0^2 U^{\rm og}_- =& 0,\quad X\in\Omega_N^-,\\
  \label{eq:uog3:n}
  u_{+}^{\rm og}(x) - U_{-}^{\rm og}(X) =& F(x), \quad x\in \Gamma,\\
  \label{eq:uog4:n}
  \partial_{\nu(x)}u_{+}^{\rm og}(x) + \gamma(x)\partial_{\nu^N(X)}U_{-}^{\rm og}(X) =& G(x),\quad x\in \Gamma,
\end{align}
where we now set $\nu^N(X)$ to be the unit normal vector of $\Gamma_N$ at $X$
pointing away from $\Omega^-_{N}$, $\gamma(x)=|M_-|^{1/2}\frac{ds(X)}{ds(x)}$,
$ds(x)$ and $ds(X)$ represent the differential arc lengths of $\Gamma$ at point
$x$ and $\Gamma_N$ at the corresponding point $X$, respectively, and
\begin{align*}
F(x)=&\lim_{Y\to X^-}U^{\rm tot}_{b,-}(Y) - \lim_{y\to x^+}u^{\rm tot}_b(y),\\
G(x)=&-\gamma(x)\lim_{Y\to X^-}\partial_{\nu^N(Y)}U^{\rm tot}_{b,-}(Y) - \lim_{y\to x^+}\partial_{\nu(x)}u^{\rm tot}_b(y).
\end{align*}
In particular, if $x$ is away from the perturbed part of $\Gamma$,
$\gamma(x)=|M_-|^{1/2}\alpha_{-}$ and (\ref{eq:utotb3}) indicates that
$F(x)=G(x)=0$ so that in fact $F$ and $G$ are compactly supported of supports
the same as the perturbed part of $\Gamma$. According to the asymptotic behavior
(\ref{eq:fp:u+}) and (\ref{eq:fp:u-}) for outgoing waves at infinity, we expect
that $u^{\rm og}_+$ and $U_-^{\rm og}$ can be absorbed by two separately defined
PMLs, as shall be discussed in the next subsection.
\subsection{Regionalized perfectly matched layers}
\begin{figure}[!htb]
\centering
\includegraphics[width=12cm]{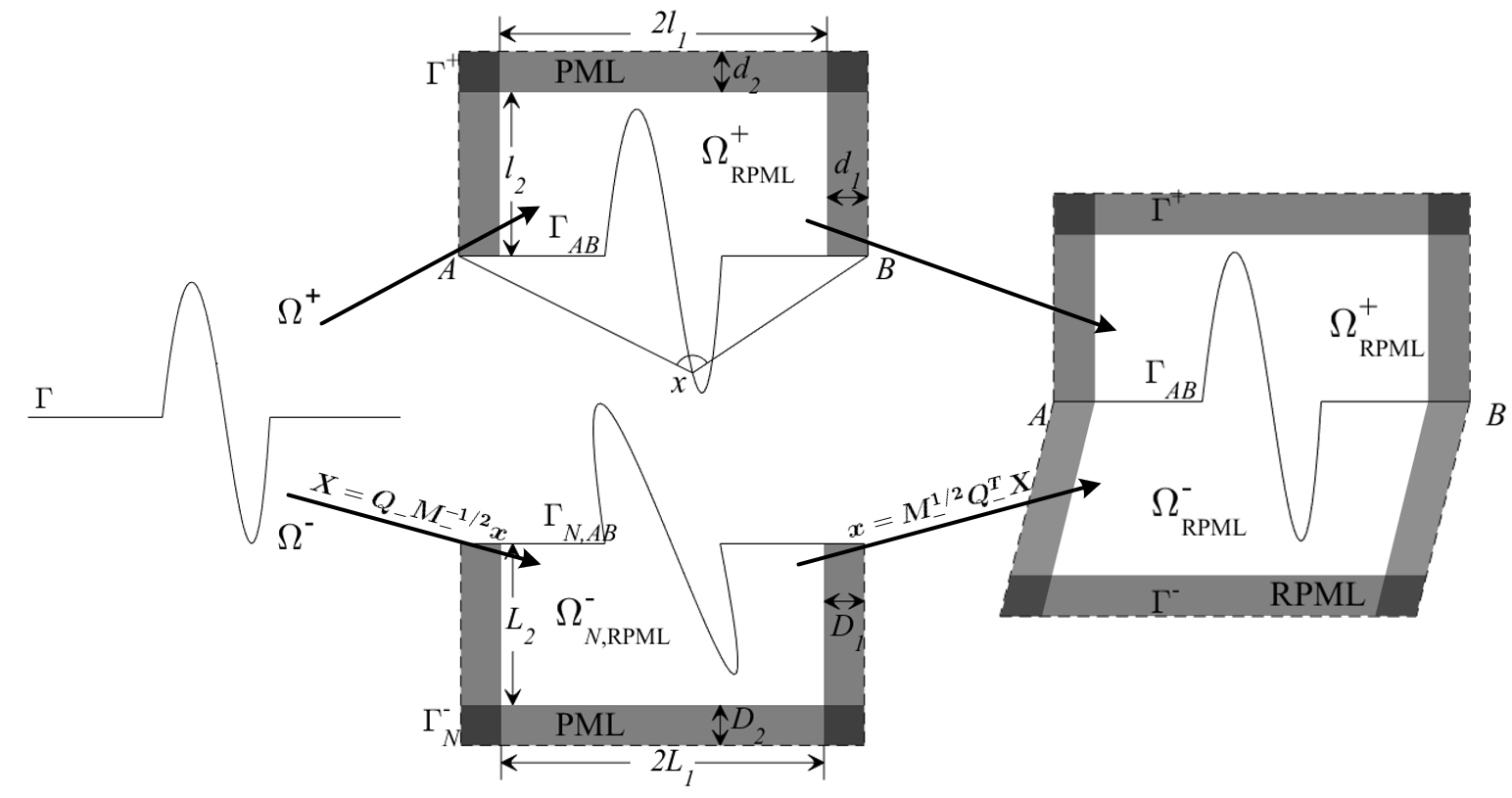}
\caption{The setup of an RPML: rectangular PMLs are setup both in $\Omega^+$ and
  in $\Omega^-_N$, the image of $\Omega^-$; $\Omega_{\rm RPML}^+$ and
  $\Omega_{N,\rm RPML}^-$ are the truncated regions; $\Omega_{\rm RPML}^-$, the
  preimage of $\Omega_{N,\rm RPML}^-$, is slanted due to the change of
  coordinates.}
\label{fig:pml}
\end{figure}

As $\Omega^+$ and $\Omega_N^-$ have different coordinate systems, we introduce
an RPML to complexify the coordinates of $x$ and $X$. As indicated by
its name, the RPML should be defined regionally and hence is no longer uniaxial.
Let $l_i$, $L_i$, $d_i$, and $D_i$ be positive constants for $i=1,2$. As shown
in Figure~\ref{fig:pml},
$\Omega^+_{\rm RPML} = \Omega^+\cap\{x:|x_1|<l_1+d_1,
|x_2|<l_2+d_2\}$, $\Omega^-_{N,\rm RPML} = \Omega^-_N\cap\{X:|X_1|<L_1+D_1,
|X_2|<L_2+D_2\}$, and $\Omega^-_{\rm RPML}$ be the preimage of $\Omega^-_{N,\rm
  RPML}$. We choose $l_i$ and $L_i$ to be sufficiently large such that
$\Omega^+_{\rm RPML}\cup\Omega^-_{\rm RPML}$ enclose the perturbed part of
$\Gamma$. Let $\Gamma^+ = \partial\Omega^+_{\rm RPML}\backslash\Gamma$,
$\Gamma_N^-= \partial\Omega^-_{N, \rm RPML}\backslash \Gamma_N$, and $\Gamma^-$
be the preimage of $\Gamma_N^-$. For $x=(x_1,x_2)^{T}$ in the top region
$\Omega^+$, we define
\begin{equation}
  \label{eq:tx}
  \tilde{x}_i = x_i + \bi\int^{x_i}_0\sigma_i(t)\mathrm{d}t,\quad i=1,2,
\end{equation}
where the two functions $\sigma_i(x)$ are positive for $|x_i|\in[l_i,l_i+d_i]$,
and are zero elsewhere. For $X=(X_1,X_2)^{T}$ in the bottom region $\Omega_N^-$,
the image of $\Omega^-$, we define
\begin{equation}
  \tilde{X}_i = X_i + \bi\int^{X_i}_0\sigma^N_i(t)\mathrm{d}t,\quad i=1,2,
\end{equation}
where the two functions $\sigma^N_i(X)$ are positive for
$|X_i|\in[L_i,L_i+d_i]$, and are zero elsewhere. The regions with nonzero
$\sigma_i$ or nonzero $\sigma_i^N$ are called the RPML regions, as indicated by
the shaded regions in Figure~\ref{fig:pml}. Thus, $d_i$ and $D_i$ represent the
thicknesses of the RPML.

Let $\tx = (\tx_1,\tx_2)^{T}$ for $x\in\Omega^+$, and $\tX=(\tX_1,\tX_2)^{T}$
for $X\in\Omega^-_N$. Analogous to Green's function in a layered isotropic
medium in \cite{chezhe10}, the background Green's function $G(x;\cdot)$ can be
analytically extended to well define $G(\tx;\cdot)$ for $x\in\Omega_{\rm
  RPML}^+$ and $G_-(\tX;\cdot)$ for $X\in\Omega_{N, \rm RPML}^-$. By the
reciprocity relation (\ref{eq:rec:rel}) and Green's representation formula
(\ref{eq:rep:GG-}), analytic continuation also applies for $u_+^{\rm og}$ and
$U_-^{\rm og}$ so that $\tu^{\rm og}_+(x)=u_+^{\rm og}(\tx)$ and $\tU_{-}^{\rm
  og}(X)=U_-^{\rm og}(\tX)$ are well-defined. By the chain rule, we obtain
\begin{align}
  \label{eq:gov:tuog}
  \nabla\cdot({\bf A}\nabla \tilde{u}_+^{og}) + k_0^2J \tu^{og}_+ &= 0,\quad {\rm in} \quad \Omega^+,\\
  \label{eq:gov:tUog}
  \nabla_X\cdot({\bf A}_N\nabla_X \tilde{U}_-^{og}) + k_0^2 J_N \tU_-^{og} &= 0,\quad {\rm in} \quad \Omega^-_N,
\end{align}
where $q_i=1+\mathbf{i}\sigma_i(x_i)$, $q^N_i=1+\mathbf{i}\sigma^N_i(X_i)$,
${\bf A}={\rm Diag}\{q_2/q_1,q_1/q_2\}$, ${\bf A}_N={\rm
  Diag}\{q_2^N/q_1^N,q_1^N/q_2^N\}$, $J=q_1(x)q_2(x)$ and $J_N=q_1^N(X)q_2^N(X)$.

Across the interface $\Gamma$ or $\Gamma_N$, we expect that the interface
conditions (\ref{eq:uog2}) and (\ref{eq:uog3}) should be continued to inside the
RPML. This relies on the following condition 
\begin{equation}
  \label{eq:cond:sigma1N}
  \sigma^N_1(X_1) = \sigma_1(x_1),\quad x_1\in\mathbb{R},
\end{equation}
as indicated by the following lemma.
\begin{mylemma}
  Under the condition (\ref{eq:cond:sigma1N}), the interface conditions
  (\ref{eq:uog3:n}) and (\ref{eq:uog4:n}) can be analytically continued to
  \begin{align}
    \label{eq:int:pml1}
    \tu_+^{\rm og}(x) - \tU_-^{\rm og}(X) =& F(x),\\ 
    \label{eq:int:pml2}
    \partial_{\nu(x)}\tu_{+}^{\rm og}(x) +\gamma(x)\partial_{\nu^N(X)}\tU_{-}^{\rm og}(X) =& G(x),
  \end{align}
 for any $x\in \Gamma$ and $X=Q_-M_-^{-1/2}x\in\Gamma_N$.
 \begin{proof}
   Outside the RPML region, (\ref{eq:int:pml1}) and (\ref{eq:int:pml2}) are exactly
   the same as (\ref{eq:uog3:n}) and (\ref{eq:uog4:n}). Inside the RPML region, $\Gamma$
   coincides with $x_2=0$ so that $X_1=\alpha_-x_1$, $\nu(x)=(0,-1)^{T}$, $\nu^N(X)=(0,1)^{T}$,
   and $F(x)=G(x)\equiv 0$. Condition (\ref{eq:cond:sigma1N}) directly implies
   that on $x_2=0$,
   \[
     \tX_1 = X_1 + \int_0^{X_1}\sigma_1^N(t)dt = \alpha_-x_1 +
     \int_{0}^{\alpha_-x_1}\sigma_1(t/\alpha_-)dt = \alpha_-\tx_1,
   \]
    so that $\tx=M_-^{1/2}Q_-^T\tX$.
   Since $F(x)=0$, the interface condition (\ref{eq:uog3:n}) reduces to
   \[
     U_-^{\rm og}(X) = u_+^{\rm og}(M_-^{1/2}Q^{T}_-X).
   \]
   The identity theorem for analytic functions directly implies 
   \[
     U_-^{\rm og}(\tX) = u_+^{\rm og}(M_-^{1/2}Q^{T}_-\tX) = u_+^{\rm og}(\tx).
   \]
   Equation (\ref{eq:int:pml2}) can be proved similarly.
 \end{proof}
\end{mylemma}
Now, we directly truncate $\tu^{\rm og}_+$ and $\tU^{\rm og}_-$ onto
$\Omega^+_{\rm RPML}$ and $\Omega^-_{N,\rm RPML}$, respectively, by imposing the
following Dirichlet boundary conditions
\begin{align}
  \label{eq:bc:uog}
  \tu^{\rm og}_+ =& 0,\quad {\rm on}\quad \Gamma^+,\\
  \label{eq:bc:Uog}
  \tU^{\rm og}_- =& 0,\quad {\rm on}\quad \Gamma_N^-.
\end{align}
Since $u^{\rm og}_+$ and $U^{\rm og}_-$ are purely outgoing waves at infinity,
the artificial boundary conditions (\ref{eq:bc:uog}) and (\ref{eq:bc:Uog}) are
expected to induce a truncation error that decays exponentially as the RPML
parameters $d_i$, $D_i$, $\sigma_i$ and $\sigma_i^N$ for $i=1,2$ increase, as
shall be validated by the numerical examples in Section 6. According to
(\ref{eq:cond:sigma1N}), we take
\begin{equation}
  \label{eq:pml:cond2}
  L_1= \alpha_-l_1,\quad D_1 = \alpha_-d_1,
\end{equation}
so that along $\Gamma$, $\tu^{\rm og}_+$ and $\tU^{\rm og}_-$ simultaneously
enter the RPML and terminate. Consequently, equations (\ref{eq:gov:tuog}),
(\ref{eq:gov:tUog}), (\ref{eq:int:pml1}), (\ref{eq:int:pml2}), (\ref{eq:bc:uog})
and (\ref{eq:bc:Uog}) form a closed boundary value problem for the two unknowns
$\tu^{\rm og}_+$ and $\tU^{\rm og}_-$. In the next subsection, we shall adopt a
previously developed BIE method to numerically solve this boundary value
problem.
\subsection{The PML-based BIE method}
As indicated by Figure~\ref{fig:pml}, let $\Gamma_{AB}$ be the truncated part of
$\Gamma$ by $\Omega^+_{\rm RPML}$ with two endpoints $A$ and $B$ and
$\Gamma_{N, AB}$ be the truncated part of $\Gamma_N$ by $\Omega_{N,\rm RPML}^-$.
It can be seen from (\ref{eq:pml:cond2}) that $\Gamma_{N,AB}$ is the image of
$\Gamma_{AB}$. We first consider $\tu^{\rm og}_+$ in $\Omega^+_{\rm RPML}$. For
any given function $g\in H^{-1/2}(\Gamma_{AB})$, consider the following boundary
value problem
\begin{align}
  \label{eq:tu+1}
  \nabla\cdot({\bf A}\nabla \tilde{u}_+^{og}) + k_0^2J \tu^{og}_+ =& 0,\quad {\rm on} \quad \Omega_{\rm RPML}^+,\\
  \tu^{\rm og}_+ =& 0,\quad {\rm on}\quad \partial\Omega^+_{\rm RPML}\backslash \Gamma,\\
  \partial_{\nu_c} \tu^{\rm og}_+|_{\Gamma_{AB}} =& g\in H^{-1/2}(\Gamma_{AB}),
\end{align}
where the conormal vector $\nu_c=\mathbf{A}^T\nu$. Fredholm theory indicates
that this problem has a unique solution $\tu^{\rm og}_+\in H^1(\Omega^+_{\rm
  RPML})$ except for $k_0$ in a countable set of eigenfrequencies;
numerically, such eigenfrequencies can be easily avoided by adjusting
$\Gamma_{AB}$ or the aforementioned RPML parameters. Then, we are able to define
a Neumann-to-Dirichlet map ${\cal N}_{AB}^+: H^{-1/2}(\Gamma_{AB})\to
\tilde{H}^{1/2}(\Gamma_{AB})$ such that $\tu^{\rm og}_+|_{\Gamma_{AB}} = {\cal
  N}_{AB}^+\partial_{\nu_c} \tu^{\rm og}_+|_{\Gamma_{AB}}$. Following
\cite{luluqia18} closely, we develop a high-accuracy BIE method to numerically
approximate ${\cal N}_{AB}^+$.

The fundamental solution of (\ref{eq:tu+1}) is 
\begin{equation}
	\label{eq:fun:helm2}
	\tilde{\Phi}(x;y) = \Phi(\tx;\ty) = \frac{\bi}{4} H_0^{(1)}[k_0\rho(\tx,\ty)],
\end{equation}
where we recall that $\Phi(x;y)=\frac{\bi}{4} H_0^{(1)}(k_0|x-y|)$,
$\ty=(\ty_1,\ty_2)$ and the complexified distance function $\rho$ is defined to
be
\begin{equation}
	\label{eq:comp:dist}
	\rho(\tilde{x},\tilde{y}) = [(\tx_1-\ty_1)^2 + (\tilde{x}_2-\tilde{y}_2)^2]^{1/2}.
\end{equation}
According to \cite{luluqia18}, we have the following Green's representation
formula
\begin{align}
	\label{eq:grerep}
	\tilde{u}_+^{\rm og}(x)=& \int_{\partial \Omega_{\rm RPML}^+}\{\tilde{\Phi}(x;y)\partial_{\nu_c}\tilde{u}_+^{\rm og}(y)-\partial_{\nu_c}\tilde{\Phi}(x;y)\tilde{u}_+^{\rm og}(y)\}ds(y)\nonumber\\
  \approx &\int_{\Gamma_{AB}}\{\tilde{\Phi}(x;y)\partial_{\nu_c}\tilde{u}_+^{\rm og}(y)-\partial_{\nu_c}\tilde{\Phi}(x;y)\tilde{u}_+^{\rm og}(y)\}ds(y),
\end{align} 
for all $x\in \Omega_1^+$. Here, we have assumed that on the RPML boundary
$\Gamma^+$, the co-normal derivative $\partial_{\nu_c}\tilde{u}_+^{\rm
  og} \approx 0$, since we expect that the RPML can absorb the outgoing wave
$\tu^{\rm og}_+$ completely. An alternative way to obtain (\ref{eq:grerep}) is
regarding itself as a direct truncation of the following integral
\begin{align}
  \label{eq:grerep2}
	\tilde{u}_+^{\rm og}(x)=& \int_{\Gamma}\{\tilde{\Phi}(x;y)\partial_{\nu_c}\tilde{u}_+^{\rm og}(y)-\partial_{\nu_c}\tilde{\Phi}(x;y)\tilde{u}_+^{\rm og}(y)\}ds(y),
\end{align} 
which does not involve the complexification of $y_2$.

As $x$ approaches $\Gamma_{AB}$, the jump relations for single- and
double- layer potentials imply \cite[Eq.~(42)]{luluqia18}
\begin{equation}
	\label{eq:pmlntd}
	{\cal K}_{AB}[\tilde{u}^{\rm og}_+](x) -{\cal K}_{0,AB} [1] (x) \tilde{u}^{\rm og}_+(x)\approx {\cal S}_{AB} [\partial_{ {\bm \nu}_c}\tilde{u}^{\rm og}_+] (x).
\end{equation}
Here, we have defined the following boundary integral operators on
$\Gamma_{AB}$,
\begin{align}
  \label{eq:def:S}
	{\cal S}_{AB}[\phi](x) &= 2\int_{\Gamma_{AB}} \tilde{\Phi}(x;y) \phi(y) ds(y),\\
  \label{eq:def:K}
	{\cal K}_{AB}[\phi](x) &= 2\,{\rm p.v.}\int_{\Gamma_{AB}} \partial_{{\nu_c}} \tilde{\Phi}(x;y) \phi(y) ds(y),\\
  \label{eq:def:K0}
	{\cal K}_{0,AB}[\phi](x) &= -\angle AxB/\pi + 2\,{\rm p.v.}\int_{\Gamma_{AB}} \partial_{{\nu_c}} \tilde{\Phi}_0(x;y) \phi(y) ds(y),
\end{align}
where the angle $\angle AxB$ is indicated in Figure~\ref{fig:pml} and
\begin{equation}
	\label{eq:green:complap}
	\tilde{\Phi}_0(x;y) = -\frac{1}{2\pi}\log\rho(\tilde{x},\tilde{y}),
\end{equation}
is the fundamental solution of the complexified Laplace equation
\begin{equation}
	\label{eq:cmp:lap}
	{\nabla\cdot({\bf A}\nabla \tilde{u})} = 0.
\end{equation}
Roughly speaking, $\tu^{\rm og}_+|_{\Gamma_{AB}} \approx ({\cal K}_{AB}-{\cal
  K}_{0,AB} [1])^{-1}{\cal S}_{AB}\partial_{\nu_c}\tu^{\rm og}_+|_{\Gamma_{AB}}$
so that the NtD operator ${\cal N}_{AB}^+\approx ({\cal K}_{AB}-{\cal K}_{0,AB}
[1])^{-1}{\cal S}_{AB}$.

Suppose the piecewise smooth curve $\Gamma_{AB}$ is parameterized by
$\{x(s)=(x_1(s), x_2(s))|\;0\leq s\leq L\}$, where $s$ is the arclength
parameter. Since corners may exist, $\tilde{u}^{\rm og}_+(x(s))$ can have corner
singularities in its derivatives at corners. To smoothen $\tilde{u}^{\rm og}_+$,
we introduce a {grading} function $s=w(t), 0\leq t\leq 1$. For a smooth segment
of $\Gamma_{AB}$ corresponding to $s\in[s^0, s^1]$ and $t\in[t^0, t^1]$ such
that $s^i=w(t^i)$ for $i=0,1$, where $s^0$ and $s^1$ correspond to two corners,
we take \cite[Eq. (3.104)]{colkre13}
\begin{equation}
	\label{eq:gfun}
	s=w(t) = \frac{s^0w_1^p + s^1w_2^p}{w_1^p+w_2^p},\quad t\in[t^0,t^1],
\end{equation}
where the positive integer $p$ ensures {that} the derivatives of $w(t)$ vanish
at the corners up to order $p$,
\[
	w_1=\left(\frac{1}{2}-\frac{1}{p}\right)\xi^3+\frac{\xi}{p}+\frac{1}{2},\quad
	w_2 = 1-w_1,\quad \xi = \frac{2 t - (t^0+t^1)}{t^1-t^0}.
\]

To simplify the {notation}, we shall use $x(t)$ to denote $x(w(t))$, and $x'(t)$
to denote $\frac{dx}{ds}(w(t))w'(t)$ in the following. Assume that {$[0,1]$} is
uniformly sampled {by $N_{\rm tot}$ grid points $\{t_j=jh\}_{j=1}^{N_{\rm tot}}$
  with even $N_{\rm tot}$ and}
grid size $h=1/N_{\rm tot}$ and that the grid points contain {all the} corner points.
Thus, ${\cal S}_{AB}[\partial_{\nu_c}\tu^{\rm og}_+]$ at point $x=x(t_j)$ can be
parameterized by
\begin{align}
	\label{eq:intS}
	{\cal S}_{AB}[\partial_{{\nu}_c}\tilde{u}^{\rm og}_+] (x(t_j)) &= \int_0^1
	S(t_j,t) \psi^{\rm s}(t) dt,
\end{align}
where $S(t_j,t)=\frac{\bi}{2} H_0^{(1)}(k\rho(x(t_j),x(t)))$, and the scaled
co-normal vector $\psi_+^{\rm s}(t)= \partial_{\nu_c}\tilde{u}^{\rm
  og}_+(x(t))|x'(t)|$, smoother than $\partial_{\nu_c}\tilde{u}^{\rm
  og}_+(x(t))$, is introduced to regularize the approximation of ${\cal
  N}_{AB}$.

Considering the logarithmic singularity of $S(t_j,t)$ at $t=t_j$, we can
discretize the integral in (\ref{eq:intS}) by Alpert's 6th-order hybrid Gauss-trapezoidal
quadrature rule \cite{alp99} and then by trigonometric interpolation to get
\begin{equation}
	{\cal S}_{AB}[\partial_{\nu_c}\tilde{u}^{\rm og}_+]\left[
		\begin{array}{c}
			x(t_1)\\
			\vdots\\
			x(t_{N_{\rm tot}})
		\end{array}
	\right] \approx {\bf S}_{AB} \left[
		\begin{array}{c}
			\psi_+^{\rm s}(t_1)\\
			\vdots\\
			\psi_+^{\rm s}(t_{N_{\rm tot}})
		\end{array}
	\right],
\end{equation}
where the $N_{\rm tot}\times N_{\rm tot}$ matrix ${\bf S}_{AB}$ approximates ${\cal S}_{AB}$. One similarly
approximates ${\cal K}_{AB}[\tu^{\rm og}_+](x(t_j))$ and ${\cal K}_{0,AB}[1](x(t_j))$ for
$j=1,\cdots, N_{\rm tot}$, so that we obtain on $\Gamma_{AB}$ that
\begin{equation}
  \label{eq:ntd:+}
  {\bm u}_+^{\rm og}\approx {\bm N}^+_{AB}{\bm \psi}_+^{\rm s},
\end{equation}
where ${\bm u}_+^{\rm og}$ and ${\bm \psi}_+^{\rm s}$ represent $N_{\rm tot}\times 1$
column vectors of $\tu_+^{\rm og}$ and $\psi_+^{\rm s}$ at the $N_{\rm tot}$ grid points
of $\Gamma_{AB}$, respectively, and ${\bm N}^+_{AB}$ represents an $N_{\rm
  tot}\times N_{\rm tot}$
matrix approximating ${\cal N}_{AB}^+$.

Now consider $\tU_-^{\rm og}$ in $\Omega_{N,\rm RPML}^-$. Let $\nu^N_c={\bf
  A}_N\nu^N$ and $\Psi_-^{\rm s}(t)=\partial_{\nu^N_c}\tU^{\rm
  og}_-(X(t))|X'(t)|$, where $X(t)$ parameterizes $\Gamma_{N,AB}$ and is the
image of $x(t)$. Following the same procedure as above, one obtains
\begin{equation}
  \label{eq:ntd:-}
  {\bm U}_-^{\rm og}\approx {\bm N}^-_{AB}{\bm \Psi}_-^{\rm s},
\end{equation}
where ${\bm U}_-^{\rm og}$ and ${\bm \Psi}_-^{\rm s}$ represent $N_{\rm tot}\times 1$
column vectors of $\tU_-^{\rm og}$ and $\Psi_-^{\rm s}$ at $N_{\rm tot}$ grid points of
$\Gamma_{N,AB}$, respectively, and ${\bm N}^-_{AB}$ is an $N_{\rm tot}\times
N_{\rm tot}$ matrix
relating to the NtD operator mapping $\partial_{\nu^N_c}\tU^{\rm og}_-$ to
$\tU^{\rm og}_-$. To match the two equations (\ref{eq:ntd:+}) and
(\ref{eq:ntd:-}) on $\Gamma_{AB}$, the $N_{\rm tot}$ grid points of $\Gamma_{N,AB}$ must
be the image of the $N_{\rm tot}$ grid points of $\Gamma_{AB}$. Consequently,
(\ref{eq:int:pml1}) and (\ref{eq:int:pml2}) imply
\begin{align}
  \label{eq:lin:sys}
  \left[
  \begin{array}{ll}
    {\bm N}_{AB}^+ & -{\bm N}_{AB}^-\\
    {\bm I}_{N_{\rm tot}} & {\bm D}_{N_{\rm tot}}\\
   \end{array}
  \right]
  \left[
  \begin{array}{l}
    {\bm \psi}_+^{\rm s}\\
    {\bm \Psi}_-^{\rm s}\\
   \end{array}
  \right] =   \left[
  \begin{array}{l}
    {\bm F}\\
    {\bm G}\\
   \end{array}
  \right],
\end{align}
where ${\bm I}_{N_{\rm tot}}$ is the $N_{\rm tot}\times N_{\rm tot}$ identity
matrix, ${\bm D}_{N_{\rm tot}}={\rm
  Diag}\{\gamma(x(t_1)),\cdots,\gamma(x(t_{N_{\rm tot}}))\}$, ${\bm F}$ and
${\bm G}$ are two $N_{\rm tot}\times 1$ vectors consist of elements $F(x(t_j))$
and $G(x(t_j))|w'(t_j)|q_1(x_1(t_j))$, respectively. Here, equation
(\ref{eq:int:pml2}) is multiplied by $|w'(t_j)|q_1(x_1(t_j))$ due to the
relation between $\psi_+^{\rm s}$ and $\partial_{\nu_c} u^{\rm og}_+$.

Solving the above linear system, we get ${\bm \psi}_+^{\rm s}$ and ${\bm
  \Psi}_-^{\rm s}$ on $\Gamma_{AB}$ and $\Gamma_{N,AB}$, respectively, and hence
obtain ${\bm u}_+^{\rm og}$ and ${\bm U}_-^{\rm og}$ by (\ref{eq:ntd:+}) and (\ref{eq:ntd:-}). Green's representation
formula (\ref{eq:grerep}) then applies to get $\tu^{\rm og}_+$ in $\Omega_{\rm
  RPML}^+$ and $\tU^{\rm og}_-$ in $\Omega^-_{N,\rm RPML}$. Consequently, $u^{\rm
  og}$ and hence $u^{\rm tot}$ become available in the physical regions of
$\Omega^+_{\rm RPML}\cup \Omega^-_{\rm RPML}$, where we recall that
$\Omega^-_{\rm RPML}$ is the preimage of $\Omega^-_{N, \rm RPML}$.
\subsection{Stability of the RPML and Uniqueness of the orthogonal matrix
  $Q_-$}
To conclude this section, \tcr{we give a formal and numerical justification
  regarding the stability of the RPML. Physically speaking, by the Green's
  representation formula (\ref{eq:rep:G}) in Theorem~\ref{thm:greenrep} and the
  outgoing behavior of $G$ indicated in Lemma~\ref{lem:asym:G}, the change of
  coordinates $X=Q_-M^{-1/2}_-x$ in (\ref{eq:qpm}) indicate that $u^{\rm
    og}_+(x)$ is purely outgoing in $\Omega^+$ in the $x$-coordinate system and
  $U^{\rm og}_-(X)$ in $\Omega_N^-$, the image of $\Omega^-$, is also purely
  outgoing in the $X$-coordinate system. It is obvious that the RPML must
  perfectly absorb $u^{\rm og}_+$ and $U^{\rm og}_-$, a rigorous proof of which
  shall be presented in \cite{lu21ort}.}

In the following, we discuss why the orthogonal matrix $Q_-$ should be chosen as
in (\ref{eq:qpm}). To illustrate this, we first give a formal proof on the
failure of the uniaxial PML (UPML), the instability of which for homogeneous and
orthotropic background has been justified in \cite{becfaujol03}. In the UPML,
the coordinate transformation (\ref{eq:tx}) applies in both $\Omega^+$ and
$\Omega^-$. Then, we can define the UPML regions $\Omega_{\rm UPML}^{\pm}=\{x:
|x_1|\leq l_1+d_1, \pm x_2<l_2+d_2\}\cap \Omega^{\pm}$. Let $\Omega_{\rm
  UPHY}^{\pm}$ and $\Omega_{\rm RPHY}^{\pm}$ be the physical regions of
$\Omega_{\rm UPML}^{\pm}$ and $\Omega_{\rm RPML}^{\pm}$, respectively. Note that
$\Omega_{\rm UPML}^+=\Omega_{\rm RPML}^+$ and $\Omega_{\rm UPHY}^+=\Omega_{\rm
  RPHY}^+$. Let $\tu^{\rm og}_{\pm, \rm UPML}$ denote the UPML-truncated wave
fields in $\Omega_{\rm UPML}^{\pm}$. Following the same procedure of the
PML-based BIE method, we obtain a linear system for unknowns
$\partial_{\nu_c}\tu^{\rm og}_{\pm, \rm UPML}$ on $\Gamma$. As the interface
conditions are only defined on $\Gamma$, such a linear system should give rise
to $\tu_{-,\rm UPML}^{\rm og}(x)=\tU_{-}^{\rm og}(X)$ and
$\partial_{\nu_c}\tu^{\rm og}_{-,\rm
  UPML}(x)=|M_-|^{1/2}\partial_{\nu^N_c}\tU^{\rm og}_-(X) ds(X)/ds(x)$ for any
$x\in \Gamma$.

Now, for any $x\in \Omega_{\rm UPHY}^-\cap \Omega_{\rm RPHY}^-$, the common
physical region of the UPML and the RPML, Green's representation formula,
analogous to (\ref{eq:grerep2}), implies 
\begin{align*}
	\tilde{u}_{-,\rm UPML}^{\rm og}(x)=& \int_{\Gamma}\{\tilde{\Phi}_{\rm UPML}(x;y)\partial_{\nu_c}\tilde{u}_{-,\rm UPML}^{\rm og}(y) -\partial_{\nu_c(y)}\tilde{\Phi}_{\rm UPML}(x;y)\tilde{u}_{-,\rm UPML}^{\rm og}(y)\}ds(y),
\end{align*}
where $\tilde{\Phi}_{\rm UPML}(x;y)=\Phi_{\rm UPML}(\tx;\ty)$ and
\[
  \Phi_{\rm UPML}(x;y) = |M_-|^{-1/2}\Phi(Q_-M_-^{-1/2}x;Q_-M_-^{-1/2}y)
\]
is the fundamental solution of (\ref{eq:Hhel}) for $M(x)\equiv M_-$. By the
change of variable $Y=Q_-M_-^{-1/2}y$,
\begin{align}
  \label{eq:UPML:u}
  \tilde{u}_{-,\rm UPML}^{\rm og}(x)=&\int_{\Gamma_N}\{\Phi(Q_-M_-^{-1/2}x;\tilde{Y}))\partial_{\nu^N_c}\tilde{U}_-^{\rm og}(Y)-\partial_{\nu^N_c(Y)}\Phi(Q_-M_-^{-1/2}x;\tilde{Y})\tilde{U}_-^{\rm og}(Y)\}ds(Y)\nonumber\\
  =&\int_{\Gamma_N}\{\Phi(X;\tilde{Y}))\partial_{\nu^N_c}\tilde{U}_-^{\rm og}(Y)-\partial_{\nu^N_c(Y)}\Phi(X;\tilde{Y})\tilde{U}_-^{\rm og}(Y)\}ds(Y).
\end{align}
On the other hand, 
\begin{align}
  \label{eq:RPML:u}
  \tU_-^{\rm og}(X)=&\int_{\Gamma_N}\{{\Phi}(\tX;\tY))\partial_{\nu^N_c}\tilde{U}_-^{\rm og}(Y)-\partial_{\nu^N_c}{\Phi}(\tX;\tY)\tilde{U}_-^{\rm og}(Y)\}ds(Y).
\end{align}
Thus, $\tilde{u}_{-,\rm UPML}^{\rm og}(x)= \tU_-^{\rm og}(X)$ since $\tX=X$. It
is surprising that the unstable UPML can unexpectedly provide accurate solutions
in part of its physical domain based on the proposed BIE method!

However, the situation changes considerably when $x\in \Omega_{\rm
  UPHY}^-\backslash\overline{\Omega_{\rm RPHY}^-}$. Equation (\ref{eq:RPML:u})
can still be analytically and stably continued to evaluate $\tU_-^{\rm og}(X)$
although $X$ now is in the RPML region of $\Omega_{\rm RPML}^-$. Consequently,
our RPML is stable in its whole physical region. On the contrary,
(\ref{eq:UPML:u}), computable though, cannot be analytically continued to the
region $\Omega_{\rm UPHY}^-\backslash\overline{\Omega_{\rm RPHY}^-}$ as the
branch cut of $\Phi(X;\tilde{Y})$, $\{X:\rho^2(X,\tilde{Y})<0\}$, is crossed
definitely. For example, if $Y_1=X_1$, $\tilde{Y}_1$ has a sufficiently large
imaginary part, and $X_2$ is sufficiently close to $Y_2$, then
\[
  \rho^2(X,\tilde{Y}) = (X_1-\tilde{Y}_1)^2 + (X_2-Y_2)^2<0.
\]
Consequently, $\tu^{\rm og}_{-,\rm UPML}(x)$ can not be an analytic function in
$\Omega_{\rm UPHY}^-$. But it contradicts our expectation that in the physical
region $x\in\Omega_{\rm UPHY}^-$, $\tu^{\rm og}_{-,\rm UPML}(x)$ coincides with
$u^{\rm og}(x)$, which must be analytic. We note that if numerical methods such
as finite difference or finite element methods are used, then much worse
numerical solutions are expected. The incorrectly enforced zero boundary
condition on the UPML boundary $\partial \Omega_{\rm UPML}^-\backslash\Gamma$
makes the truncation error propagate back to the whole computational domain!

In a similar fashion to the above, it can be seen that choosing $Q_-$ other than
(\ref{eq:qpm}) does not work either. Different from the UPML, a different choice
of $Q_-$ leads to a rotation of the RPML region of $\Omega_{\rm RPML}^-$. Again,
branch cut of $\Phi$ appears inside the corresponding physical region, making
such an RPML not work either.

To illustrate the above statements more clearly, we consider a specific example
below. Let $\Gamma=\{x:x_2=0\}$, $k_0=2\pi$, $\epsilon^+$ be the $3\times 3$
identity matrix $I_3$, and
\begin{equation*}
\epsilon^-=\begin{bmatrix}

     4 & 3 & 0 \\

     3 & 4 & 0 \\
     
     0 & 0 & 1 \\

	\end{bmatrix}.
\end{equation*}
We compute $u^{\rm tot}(x;x^*)$ for a cylindrical incident wave $u^{\rm
  inc}(x;x^*)$ excited by a source at $x^*=(0,0.1)^{T}$ based on three different
types of PMLs: (1) UPML; (2) RPML-I with $Q_-=I_2$; (3) RPML-II with $Q_-$
defined in (\ref{eq:qpm}). In $\Omega^+$, we let $l_1=1$ and $d_1=1$ so that
$\Gamma_{AB}=\{x:|x_1|<2, x_2=0\}$, and choose (\ref{eq:sigma1}) with $S=2$ to
complexify $x_1$ in $\Omega^+$ and $X_1$ in $\Omega_N^-$. The computation domain
is set to be $\Omega_{D}=(-2,2)\times(-2,2)$, and we let $l_2$ and $L_2$
sufficiently large so that $\Omega_D$ contains only regions that complexify
$x_1$. To obtain sufficiently accurate numerical solutions, we discretize
$\Gamma_{AB}$ by $N_{\rm tot}=800$ points.

\begin{figure}[!htb]
\centering
a)\includegraphics[width=3.25cm]{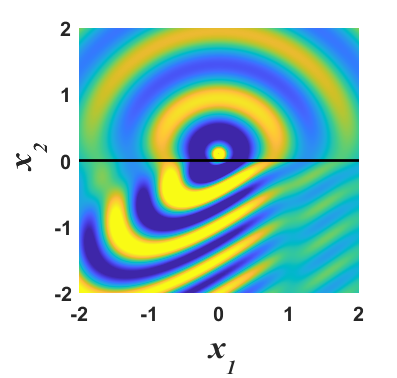}
b)\includegraphics[width=3.25cm]{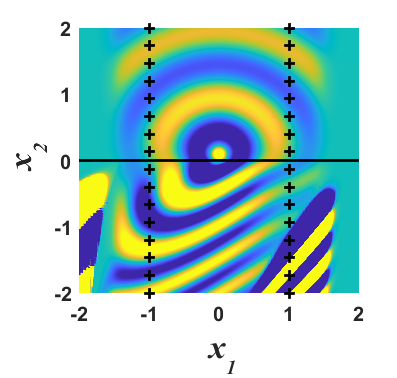}
c)\includegraphics[width=3.25cm]{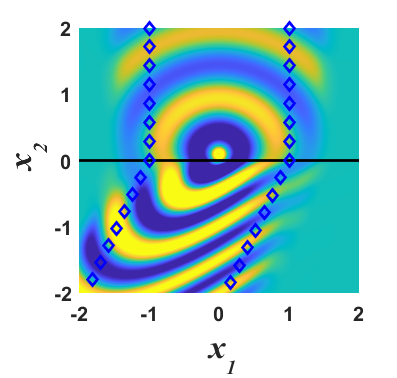}
d)\includegraphics[width=3.25cm]{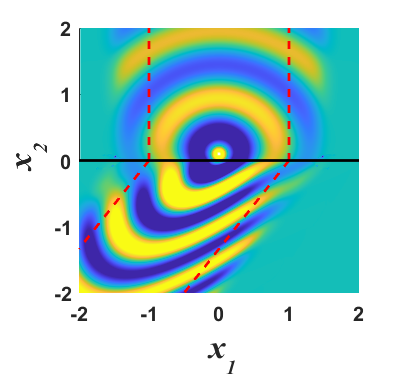}
\caption{Numerical solutions by different PMLs: a) Exact solution $u^{\rm
    tot}_{\rm exa}$; b) UPML solution $\tu_{\rm UPML}^{\rm tot}$ with PML
  entrances marked by '+'s; c) RPML-I solution $\tu_{\rm RPML,1}^{\rm tot}$ with
  entrances marked by '$\diamond$'s; d) RPML-II solution $\tu_{\rm RPML,2}^{\rm
    tot}$ with entrances marked by dashed lines.}
\label{pic:comp}
\end{figure}

\begin{figure}[!htb]
\centering
a)\includegraphics[width=4.0cm]{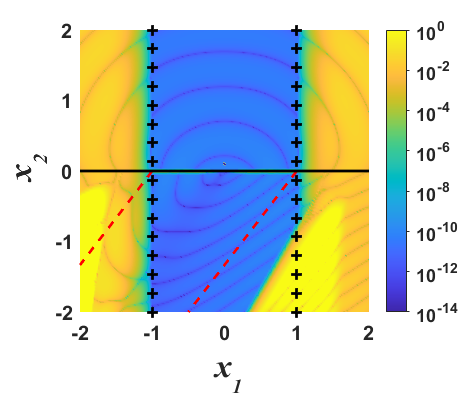}
b)\includegraphics[width=4.0cm]{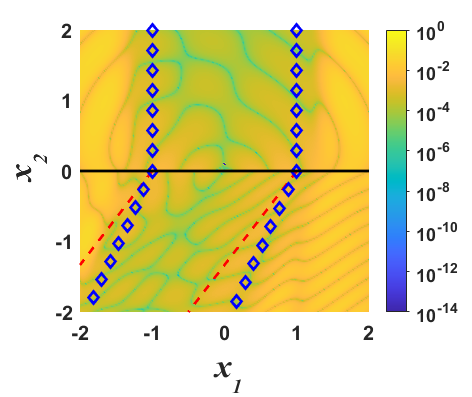}
c)\includegraphics[width=4.0cm]{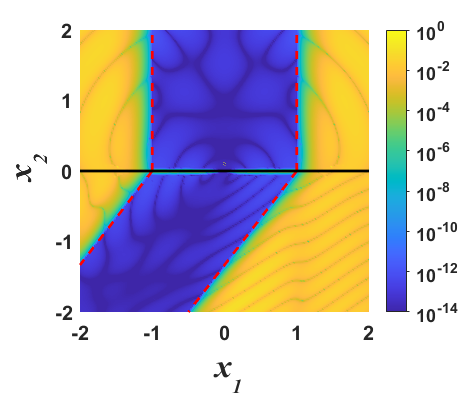}
\caption{Errors of the three numerical solutions: a) $\tu^{\rm tot}_{\rm UPML}$;
  b) $\tu^{\rm tot}_{\rm RPML,1}$; c) $\tu^{\rm tot}_{\rm RPML,2}$. Dashed lines
  in all figures indicate the entrance of the RPML-II. The '+'s indicate
  the entrance of the UPML, whereas the '$\diamond$'s indicate the entrance of
  the RPML-I.}
\label{pic:comp:error}
\end{figure}

Clearly, $u^{\rm tot}_{\rm exa}(x;x^*)=G(x;x^*)$ defined in (\ref{eq:G:++}) and
(\ref{eq:G:+-}) is the exact solution. Let $\tu^{\rm tot}_{\rm UPML}$, $\tu^{\rm
  tot}_{\rm RPML,1}$, and $\tu^{\rm tot}_{\rm RPML,2}$ be the three numerical
solutions produced by UPML, RPML-I and RPML-II. Real parts of the exact solution
and the three numerical solutions are shown in Figure~\ref{pic:comp}.
In comparison with the exact solution, the errors of the three numerical
solutions are depicted in Figure~\ref{pic:comp:error}.
As can be seen, only $\tu^{\rm tot}_{\rm RPML,2}$ perfectly coincides with the
exact solution $u^{\rm tot}_{\rm exa}$ in its physical region with pointwise
errors around ${\cal O}(10^{-12})$. In contrast, $\tu^{\rm tot}_{\rm UPML}$ and
$\tu^{\rm tot}_{\rm RPML,1}$ are accurate only when $x$ lies in their physical
regions and also in the physical region of RPML-II, as expected according to the
above theory.
\section{Numerical examples}
In this section, we carry out several numerical experiments to validate the
stability of the RPML and to illustrate the high accuracy of the proposed
BIE method. In all examples, we let the freespace wavelength $\lambda=1$ so that $k_0=2\pi$,
$\epsilon^+=I_3$, and
\begin{equation*}
\epsilon^-=\begin{bmatrix}

     4 & 1 & 0 \\

     1 & 9 & 0 \\
     
     0 & 0 & 1 \\

	\end{bmatrix}.
\end{equation*}
In the setup of the RPML, we choose 
\begin{equation}
	\label{eq:sigma1}
	\sigma_1^N(X_1)=\sigma_1(x_1) = 
	\left\{
		\begin{array}{lc}
			\frac{2Sf_{1}^{6}}{f_1^{6}+f_2^{6}}, & l_1\leq x_1\leq l_1+d_1;\\
			\sigma_1(-x_1), & -l_1-d_1\leq x_1 \leq l_1;\\
      0, & {\rm elsewhere},
		\end{array}
	\right.  
\end{equation}
where
\[
	f_1=\left(\frac{1}{2}-\frac{1}{p}\right)\bar{x}_1^3+\frac{\bar{x}_1}{p}+\frac{1}{2},\quad
	f_2 = 1-f_1,\quad \bar{x}_1 = \frac{ x_1 - (l_1+d_1)}{d_1},
\]
and $S>0$ determines the RPML strength for absorbing outgoing waves. The
function $\sigma_1$ is of class $C^6$ at $x_1=\pm l_1$, so that $(\pm
l_1,0)^{T}$ can be considered as smooth points of $\Gamma_{AB}$. We point out
that $\sigma_2$ and $\sigma_{2}^N$ are useless in the computations. We choose
$p=6$ to define the grading function $s=w(t)$ in (\ref{eq:gfun}). In each
example, to quantify the truncation error due to the RPML, we evaluate the
relative error
\[
  E_{\rm rel} = \frac{||{\bm u}_{\rm num}^{\rm tot} - {\bm u}^{\rm tot}_{\rm
      exa}||_{\infty}}{||{\bm u}_{\rm exa}^{\rm tot}||_{\infty}},
\]
where ${\bm u}_{\rm num}^{\rm tot}$ denotes the vector of a numerical solution
$u_{\rm num}^{\rm tot}$ for the total field $u^{\rm tot}$ at grid points of the
perturbed part of $\Gamma$, and ${\bm u}^{\rm tot}_{\rm exa}$ denotes the vector
of a reference solution $u_{\rm exa}^{\rm tot}$, the exact solution of the total
field $u^{\rm tot}$ if available or a sufficiently accurate numerical solution,
at the same grid points.

\textbf{Example 1.} In this example, we assume again that $\Gamma=\{x:x_2=0\}$
to check the performance of our RPML and the convergence order of the PML-based
BIE method. Here, $l_1=1$ and the PML thickness in $\Omega^+$ is fixed as
$d_1=1$ so that $\Gamma_{AB}=\{x\in\Gamma:|x_1|<2\}$ with $A=(-2,0)^{T}$ and
$B=(2,0)^{T}$.

\begin{figure}[!htb]
\centering
a)\includegraphics[width=3.3cm]{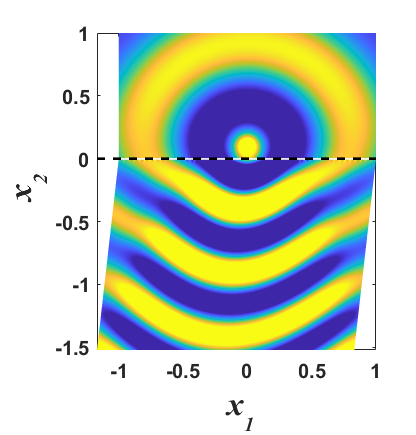}
b)\includegraphics[width=3.3cm]{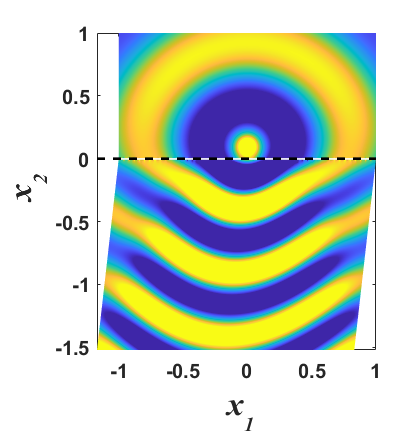}
c)\includegraphics[width=3.15cm]{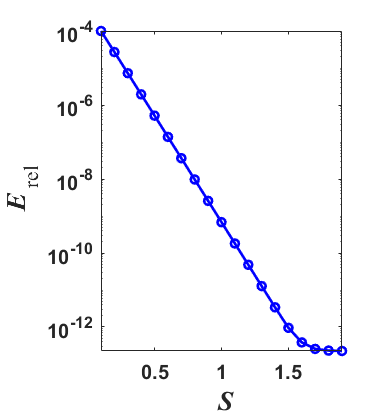}
d)\includegraphics[width=3.15cm]{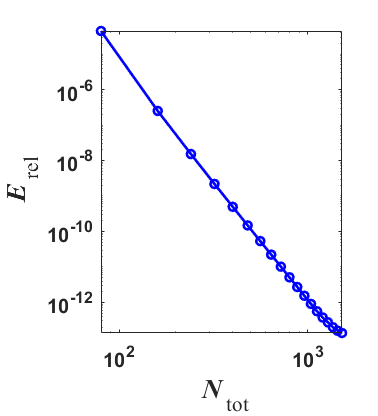}
\caption{Example 1: a) The exact solution $u^{\rm tot}_{\rm exa}$; b) An
  accurate numerical solution $u^{\rm tot}_{\rm num}$; c) $E_{\rm rel}$ against
  $S$; d) $E_{\rm rel}$ against $N_{\rm tot}$. Dashed lines in a) and b)
  indicate the interface $\Gamma$.}
\label{pic:exact}
\end{figure}

We consider only a cylindrical incident wave excited by a source located at $x^*
= (0,0.1)^{T}$. The exact solution is the Green's function $G(x;x^*)$ given by
(\ref{eq:G:++}) and (\ref{eq:G:+-}). The ``perturbed'' part on $\Gamma$ is
assumed to be $\{x:|x_1|<1,x_2=0\}$, i.e., the physical part of $\Gamma_{AB}$.
As shown in Figure \ref{pic:exact}, the computational domain $\Omega_{\rm PHY}$
is set to be the union of a rectangular region above $x_2=0$, and a slanted
region below $x_2=0$ since the physical region of $\Omega_{\rm RPML}^-$ is
always slanted. By choosing $N_{\rm tot}=1600$ grid points on $\Gamma_{AB}$, we
compute $u^{\rm tot}_{\rm num}$ on $\Omega_{\rm PHY}$, and compare it with the
exact solution $u_{\rm exa}^{\rm tot}$. Figure \ref{pic:exact} (a) and (b)
show the real parts of $u^{\rm tot}_{\rm exa}$ and $u^{\rm tot}_{\rm num}$ in
$\Omega_{\rm PHY}$. It can be seen that the two are distinguishable.

To show the stability of the RPML, we fix $N_{\rm tot}=1600$ and compute the
relative error $E_{\rm rel}$ for $S$ ranging from $0.1$ to $2$. The results are
shown in Figure \ref{pic:exact}(c), where only the vertical axis is
logarithmically scaled. We observe that $E_{\rm rel}$ decays exponentially at
the beginning and then yields to the discretization error which dominates the
relative error for large $S$. Next, we study $E_{\rm rel}$ against $N_{\rm tot}$
for $S=2$. For $N_{\rm tot}$ ranging from 80 to 1520 of step size 80, the
relative errors are depicted in Figure \ref{pic:exact}(d) where both axes are
logarithmically scaled. The slope of the decreasing part of the curve reveals
that the convergence order of the PML-based BIE method is approximately seven.
We observe that the numerical solutions are accurate to at least $12$
significant digits.

\begin{figure}[!htb]
\centering
a)\includegraphics[width=3.3cm]{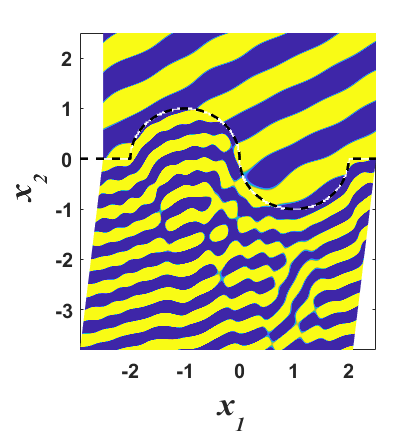}
b)\includegraphics[width=3.3cm]{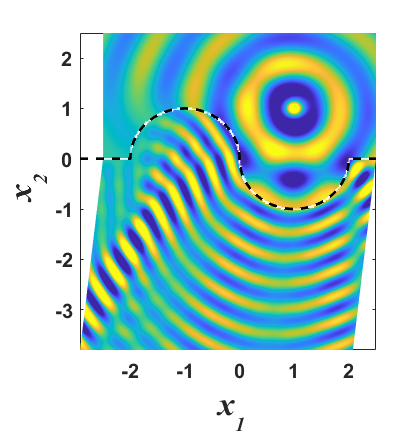}
c)\includegraphics[width=3.15cm]{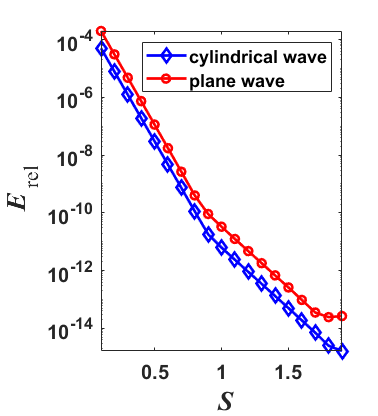}
d)\includegraphics[width=3.15cm]{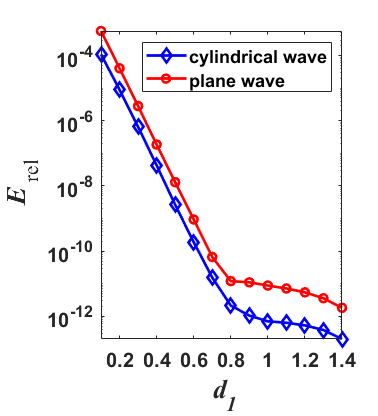}
\caption{Example 2: a) Reference solution for a plane incident wave of angle
  $\theta = \pi/3$; b) Reference solution for a cylindrical incident wave
  excited by a source at $x^*=(1,1)^{T}$; c) $E_{\rm rel}$ against $S$; d)
  $E_{\rm rel}$ against $d_1$. Dashed lines in a) and b) indicate the interface
  $\Gamma$.}
\label{pic:circle}
\end{figure}

\textbf{Example 2.} In this example, we assume that the perturbed part of
$\Gamma$ consists of two connected semicircles of radius $1$, as shown by the
dashed lines in Figure \ref{pic:circle}.

We consider two types of incidences, a plane incident wave $u^{\rm
  inc}(x;\theta)=e^{\bi k_0(\cos\theta x_1-\sin\theta x_2)}$ with
$\theta=\frac{\pi}{3}$ and a cylindrical incident wave excited by a source at
the point $x^*=(1,1)^{T}$. Here, $\Gamma_{AB}$ consists of four smooth segments. We
choose $S = 2$, $l_1=1.5$, $d_1=1.5$, and 800 points on each smooth segment so
that $N_{\rm tot}=3200$, to compute a reference solution $u^{\rm tot}_{\rm exa}$ for
either of the two incidences. Real parts of the two reference solutions are
shown in Figure \ref{pic:circle} (a) and (b). Note that the computational domain
is still slanted below $x_2=0$.

With the reference solutions available for the two incident waves, we show the
stability of the RPML. We fix $N_{\rm tot}=3200$, and check the relation between
the relative error $E_{\rm rel}$ and one of the two RPML parameters $d_1$ and
$S$ with the other one fixed. For $d_1=1.5$, the relative error $E_{\rm rel}$
against $S$, ranging from $0.1$ to $1.5$, is shown in
Figure~\ref{pic:circle}(c). For $S=2$, the relative error $E_{\rm rel}$ against
$d_1$, ranging from $0.1$ to $1.2$, is shown in Figure~\ref{pic:circle}(d). The
vertical axes in both figures are logarithmically scaled. We observe that
$E_{\rm rel}$ decays exponentially as either $S$ or $d_1$ increases for both the
two incident waves. The convergence curves indicate that the numerical solutions
are accurate to at least $11$ significant digits.

\begin{figure}[!htb]
\centering
a)\includegraphics[width=3.3cm]{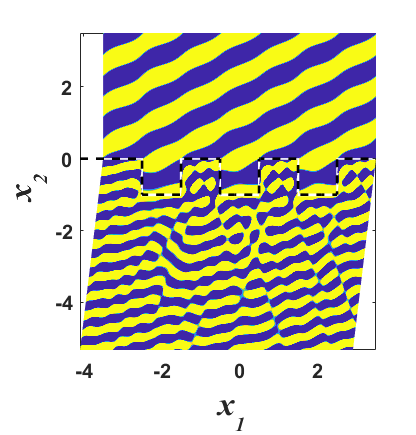}
b)\includegraphics[width=3.3cm]{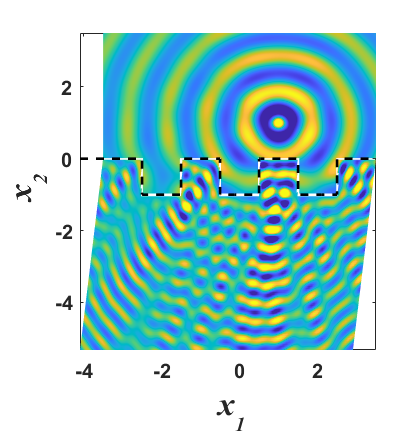}
c)\includegraphics[width=3.15cm]{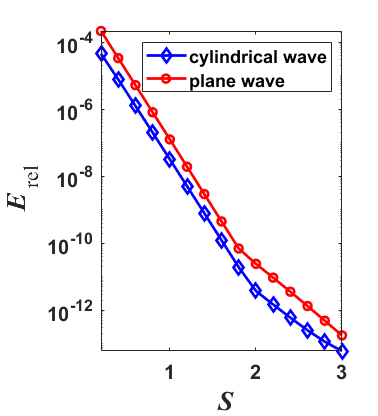}
d)\includegraphics[width=3.15cm]{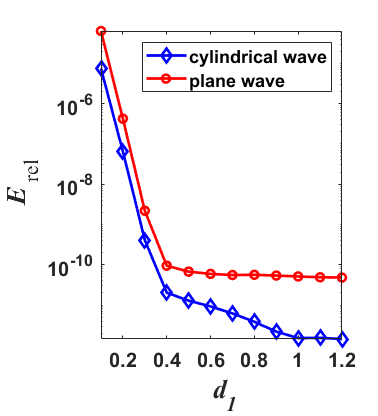}
\caption{Example 3: a) Reference solution for a plane incident wave of incident
  angle $\theta = \pi/3$; b) Reference solution for a cylindrical incident wave
  excited by $x^*=(0,1)^{T}$; c) $E_{\rm rel}$ against $S$; d) $E_{\rm rel}$
  against $d_1$. Dashed lines in a) and b) indicate the interface $\Gamma$..}
\label{pic:step}
\end{figure}

\textbf{Example 3.} In the last example, we assume that $\Omega^+$ contains
three indentations, each of which is a square of size $1$, as shown by the
dashed line in Figure \ref{pic:step}.

We consider two types of incidences, a plane incident wave $u^{\rm
  inc}(x;\theta)=e^{\bi k_0(\cos\theta x_1-\sin\theta x_2)}$ with
$\theta=\frac{\pi}{3}$ and a cylindrical incident wave excited by a source at
$x^*=(0,1)^{T}$. We choose $S = 4$, $l_1=3.5$, $d_1=1.5$, and 200 points on each of
the 13 smooth segments of $\Gamma_{AB}$ so that $N_{\rm tot}=2600$, to compute a reference
solution $u^{\rm tot}_{\rm exa}$ for either of the two incidences. Real parts of
the two reference solutions are shown in Figure \ref{pic:step} (a) and (b).

With the reference solutions available, we show the stability of the RPML now. We
fix $N_{\rm tot}=2600$, and check the relation between the relative error
$E_{\rm rel}$ and one of the two RPML parameters $d_1$ and $S$. For $d_1=1.5$,
the relative error $E_{\rm rel}$ against $S$, ranging from $0.2$ to $3.0$, is
shown in Figure~\ref{pic:circle}(c). For $S=4$, the relative error $E_{\rm rel}$
against $d_1$, ranging from $0.1$ to $1.2$, is shown in
Figure~\ref{pic:circle}(d). We observe from the two figures that $E_{\rm rel}$
decays exponentially as either $S$ or $d_1$ increases for both the two incident
waves. The convergence curves indicate that the numerical solutions are accurate
to at least $10$ significant digits. The solutions are less accurate compared
with the previous examples, since $\Gamma$ contains more corners and each smooth
segment is discretized by much less number of grid points.
 \section{Conclusion}

 In this paper, we studied wave scattering in a two-layer orthotropic medium in
 two dimensions. A novel SRC condition was proposed and a stable RPML technique
 was developed to truncate the unbounded domain. The resulting boundary value
 problem was solved by a recently developed PML-based BIE method
 \cite{luluqia18}. Numerical experiments have justified the accuracy of the
 numerical method and the stability of the RPML method, showing that the
 truncation error due to the RPML decays exponentially as the RPML
 parameters increase.

 As we can see from the numerical results, due to the transition matrix
 $Q_-M_-^{-1/2}$, the physical region of the computational domain always
 contains a slanted region. If the slanted region is too narrow, then the
 resulting computational domain will be too small, making numerical solutions
 probably useless in practice. A possible remedy could be using more generalized
 complexifications of $x_1$ and $X_1$ to enlarge the physical domain
 \cite{bonchaflihazpertja21}. We shall investigate this issue in a future work.
 Besides, we shall rigorously prove the well-posedness of the scattering problem
 and shall justify the exponentially decaying truncation error due to the RPML
 in the subsequent work \cite{lu21ort}.

 \tcr{It can be seen that the setup of our RPML does not depend on the
   wavenumber $k_0$, so that its extension to time domain is straightforward.}
 Moreover, our RPML technique exhibits deep potential in terminating waves in
 more complicated anisotropic backgrounds. Thus, it is of great interests to
 investigate the extension of the RPML to more general anisotropic media for
 both EM and elastic waves in the future.

\section*{Acknowledgements} 
W. L. would like to express his sincere gratitude to Prof. Anne-Sophie
Bonnet-BenDhia for sharing her slides, originally presented in the conference of
WAVES 2019 in Vienna, which greatly inspire the current work.
\bibliographystyle{plain}
\bibliography{wt}
\end{document}